\documentclass[12pt]{amsart}       %was 12pt
\usepackage{amssymb}
\usepackage{eucal}
\usepackage{graphicx}
\usepackage{amsmath}
\usepackage{amscd}
\usepackage[all]{xy}           %xypic macro for latex2.09

\usepackage{amsfonts,latexsym}
\usepackage{xspace}
\usepackage{hyperref}
\usepackage{epsfig}
\usepackage{float}
\usepackage{axodraw}
\usepackage{color}
\usepackage{fancybox}
\usepackage{colordvi}
\usepackage{multicol}
\usepackage{colordvi}
\newcommand{\nc}{\newcommand}
\newcommand{\delete}[1]{}
\newcommand{\deleted}[1]{}
\nc{\dfootnote}[1]{{}}          %{{}}
\nc{\ffootnote}[1]{\dfootnote{#1}}
\nc{\mfootnote}[1]{\footnote{#1}} % Use this to show footnotes
\nc{\ofootnote}[1]{\footnote{\tiny Older version: #1}} % Use this to show footnotes

\nc{\mlabel}[1]{\label{#1}}  % Use this to suppress names
\nc{\mcite}[1]{\cite{#1}}  % Use this to suppress names
\nc{\mref}[1]{\ref{#1}}  % Use this to suppress names
\nc{\mbibitem}[1]{\bibitem{#1}} % Use this to show number name

\delete{
\nc{\mlabel}[1]{\label{#1}  % Use the next two lines to show names
{\hfill \hspace{1cm}{\bf{{\ }\hfill(#1)}}}}
\nc{\mcite}[1]{\cite{#1}{{\bf{{\ }(#1)}}}}  % Use this lines to show names
\nc{\mref}[1]{\ref{#1}{{\bf{{\ }(#1)}}}}  % Use this lines to show names
\nc{\mbibitem}[1]{\bibitem[\bf #1]{#1}} % Use this to show name
}
\nc{\mkeep}[1]{\marginpar{{\bf #1}}} % Use this to show marginpar
\newtheorem{theorem}{Theorem}[section]
\newtheorem{prop}[theorem]{Proposition}
\newtheorem{defn}[theorem]{Definition}
\newtheorem{lemma}[theorem]{Lemma}
\newtheorem{coro}[theorem]{Corollary}
\newtheorem{prop-def}{Proposition-Definition}[section]

\newtheorem{tempex}[theorem]{Example}
\newtheorem{tempexs}[theorem]{Examples}
\newtheorem{temprmk}[theorem]{Remark}
\newtheorem{tempexer}{Exercise}[section]
\newenvironment{exam}{\begin{tempex}\rm}{\end{tempex}}

\nc{\vsa}{\vspace{-.1cm}} \nc{\vsb}{\vspace{-.2cm}}
\nc{\vsc}{\vspace{-.3cm}} \nc{\vsd}{\vspace{-.4cm}}
\nc{\vse}{\vspace{-.5cm}}

\delete{

}

\def\tb2{\includegraphics[scale=0.42]{tree2}}
\def\tc3{\includegraphics[scale=0.42]{tree3}}

\def\ta1{\includegraphics[scale=0.42]{tree1}}
\def\tb2{\includegraphics[scale=0.42]{tree2}}
\def\tc3{\includegraphics[scale=0.42]{tree3}}
\def\td31{\!\!\includegraphics[scale=0.42]{tree31}}
\def\te4{\includegraphics[scale=0.42]{tree4}}
\def\tf41{\!\!\includegraphics[scale=0.42]{tree41}}
\def\tg42{\!\!\includegraphics[scale=0.42]{tree42}}
\def\th43{\!\!\includegraphics[scale=0.42]{tree43}}
\def\ti5{\includegraphics[scale=0.42]{tree5}}
\def\tj51{\!\!\includegraphics[scale=0.42]{tree51}}
\def\tk52{\!\!\includegraphics[scale=0.42]{tree52}}
\def\tl53{\!\!\includegraphics[scale=0.42]{tree53}}
\def\tm54{\!\!\includegraphics[scale=0.42]{tree54}}
\def\tn55{\!\!\includegraphics[scale=0.42]{tree55}}
\def\tp56{\!\!\includegraphics[scale=0.42]{tree56}}
\def\tq57{\!\!\includegraphics[scale=0.42]{tree57}}
\def\tr58{\!\!\includegraphics[scale=0.42]{tree58}}
\def\graph1{\includegraphics[scale=0.42]{graph1}}
\def\sta1{\includegraphics[scale=0.62]{tree1}}
\def\stb2{\includegraphics[scale=0.62]{tree2}}
\def\stc3{\includegraphics[scale=0.62]{tree3}}
\def\std31{\!\!\includegraphics[scale=0.62]{tree31}}
\def\ste4{\includegraphics[scale=0.62]{tree4}}
\def\stf41{\!\!\includegraphics[scale=0.62]{tree41}}
\def\stg42{\!\!\includegraphics[scale=0.62]{tree42}}
\def\sth43{\!\!\includegraphics[scale=0.62]{tree43}}
\def\sti5{\includegraphics[scale=0.62]{tree5}}
\def\stj51{\!\!\includegraphics[scale=0.62]{tree51}}
\def\stk52{\!\!\includegraphics[scale=0.62]{tree52}}
\def\stl53{\!\!\includegraphics[scale=0.62]{tree53}}
\def\stm54{\!\!\includegraphics[scale=0.62]{tree54}}
\def\stn55{\!\!\includegraphics[scale=0.62]{tree55}}
\def\stp56{\!\!\includegraphics[scale=0.62]{tree56}}
\def\stq57{\!\!\includegraphics[scale=0.62]{tree57}}
\def\str58{\!\!\includegraphics[scale=0.62]{tree58}}
\delete{

}
\def\coprod2{\includegraphics[scale=0.35]{coprod2}}
\def\l2cop{\includegraphics[scale=0.3]{coprod2l}}
\def\c2cop{\includegraphics[scale=0.35]{coprod2c}}
\def\r2cop{\includegraphics[scale=0.3]{coprod2r}}

\def\lcop1{\includegraphics[scale=0.5]{coprod1l}}
\def\rcop1{\includegraphics[scale=0.5]{coprod1r}}
\def\rr1cop{\includegraphics[scale=0.5]{coprod1rr}}
\def\rl1cop{\includegraphics[scale=0.5]{coprod1rl}}

\def\sshu1{\!\!\includegraphics[scale=0.41]{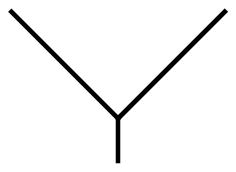}}
\def\sa2tree{\!\!\includegraphics[scale=0.41]{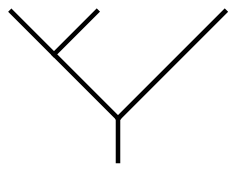}}
\def\sb3tree{\!\!\includegraphics[scale=0.41]{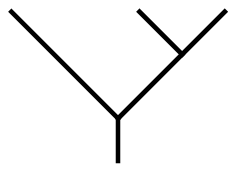}}
\def\sc4tree{\!\!\includegraphics[scale=0.41]{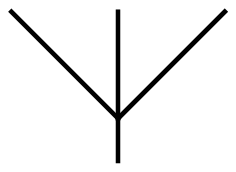}}

\def\shu1{\!\!\includegraphics[scale=0.51]{1tree.eps}}
\def\shub2{\!\!\includegraphics[scale=0.51]{2tree.eps}}
\def\shuc3{\!\!\includegraphics[scale=0.51]{3tree.eps}}
\def\shud4{\!\!\includegraphics[scale=0.51]{4tree.eps}}
\def\5tree{\!\!\includegraphics[scale=0.51]{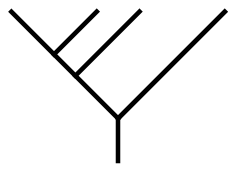}}
\def\6tree{\!\!\includegraphics[scale=0.51]{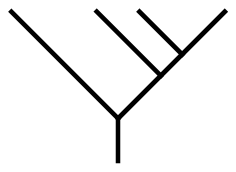}}
\def\7tree{\!\!\includegraphics[scale=0.51]{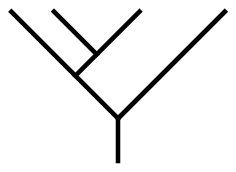}}
\def\8tree{\!\!\includegraphics[scale=0.51]{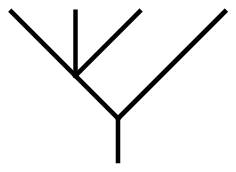}}
\def\9tree{\!\!\includegraphics[scale=0.51]{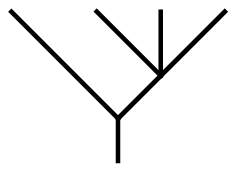}}
\def\a1tree{\!\!\includegraphics[scale=0.51]{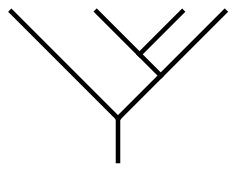}}
\def\b1tree{\!\!\includegraphics[scale=0.51]{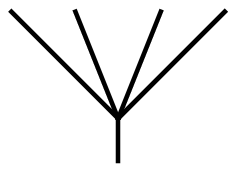}}
\def\c1tree{\!\!\includegraphics[scale=0.51]{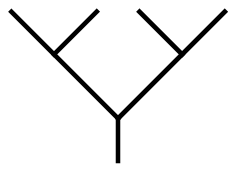}}
\def\d1tree{\!\!\includegraphics[scale=0.51]{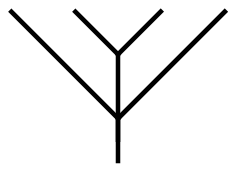}}
\def\e1tree{\!\!\includegraphics[scale=0.51]{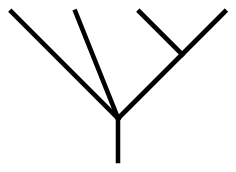}}
\def\f1tree{\!\!\includegraphics[scale=0.51]{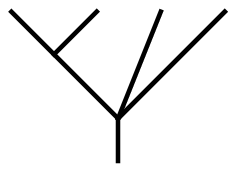}}
\def\ta1{\includegraphics[scale=0.42]{tree1}}
\def\tb2{\includegraphics[scale=0.42]{tree2}}
\def\tc3{\includegraphics[scale=0.42]{tree3}}
\def\td31{\!\!\includegraphics[scale=0.42]{tree31}}
\def\te4{\includegraphics[scale=0.42]{tree4}}
\def\tf41{\!\!\includegraphics[scale=0.42]{tree41}}
\def\tg42{\!\!\includegraphics[scale=0.42]{tree42}}
\def\th43{\!\!\includegraphics[scale=0.42]{tree43}}
\def\ti5{\includegraphics[scale=0.42]{tree5}}
\def\tj51{\!\!\includegraphics[scale=0.42]{tree51}}
\def\tk52{\!\!\includegraphics[scale=0.42]{tree52}}
\def\tl53{\!\!\includegraphics[scale=0.42]{tree53}}
\def\tm54{\!\!\includegraphics[scale=0.42]{tree54}}
\def\tn55{\!\!\includegraphics[scale=0.42]{tree55}}
\def\tp56{\!\!\includegraphics[scale=0.42]{tree56}}
\def\tq57{\!\!\includegraphics[scale=0.42]{tree57}}
\def\tr58{\!\!\includegraphics[scale=0.42]{tree58}}

\def\dec1tree{\!\!\includegraphics[scale=0.51]{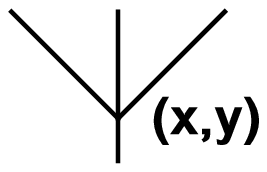}}
\def\xtree{\!\!\includegraphics[scale=0.41]{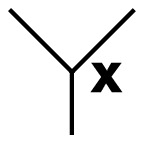}}
\def\xyztree{\!\!\includegraphics[scale=0.41]{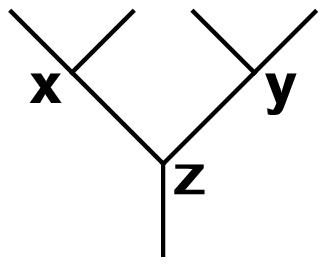}}

\def\ta1{\includegraphics[scale=0.42]{tree1}}
\def\tb2{\includegraphics[scale=0.42]{tree2}}
\def\tc3{\includegraphics[scale=0.42]{tree3}}
\def\td31{\!\!\includegraphics[scale=0.42]{tree31}}
\def\te4{\includegraphics[scale=0.42]{tree4}}
\def\tf41{\!\!\includegraphics[scale=0.42]{tree41}}
\def\tg42{\!\!\includegraphics[scale=0.42]{tree42}}
\def\th43{\!\!\includegraphics[scale=0.42]{tree43}}
\def\ti5{\includegraphics[scale=0.42]{tree5}}
\def\tj51{\!\!\includegraphics[scale=0.42]{tree51}}
\def\tk52{\!\!\includegraphics[scale=0.42]{tree52}}
\def\tl53{\!\!\includegraphics[scale=0.42]{tree53}}
\def\tm54{\!\!\includegraphics[scale=0.42]{tree54}}
\def\tn55{\!\!\includegraphics[scale=0.42]{tree55}}
\def\tp56{\!\!\includegraphics[scale=0.42]{tree56}}
\def\tq57{\!\!\includegraphics[scale=0.42]{tree57}}
\def\tr58{\!\!\includegraphics[scale=0.42]{tree58}}

\def\tde31{\!\!\includegraphics[scale=0.42]{tree31e.eps}}
\def\tdn31{\!\!\includegraphics[scale=0.42]{tree31n.eps}}

\def\ptf41{\!\!\includegraphics[scale=0.42]{planetreet41}}
\def\ptj51{\!\!\includegraphics[scale=0.42]{planetreet51}}
\def\ptk52{\!\!\includegraphics[scale=0.42]{planetreet52}}
\def\ptp56{\!\!\includegraphics[scale=0.42]{planetreet53}}

\def\dtd31{\!\!\includegraphics[scale=0.42]{dtd31}}
\def\ldtd31{\!\!\includegraphics[scale=0.42]{ldtd31}}

\def\ltde31{\!\!\includegraphics[scale=0.5]{ltree31e.eps}}
\def\ltdn31{\!\!\includegraphics[scale=0.5]{ltree31n.eps}}
\def\ltf41{\includegraphics[scale=0.5]{l41.eps}}
\def\lti4{\includegraphics[scale=0.5]{l411.eps}}
\def\ltg42{\includegraphics[scale=0.5]{l42.eps}}
\def\lth43{\includegraphics[scale=0.5]{l43.eps}}

\def\lta1{\includegraphics[scale=0.5]{ltree1.eps}}
\def\ltb2{\includegraphics[scale=0.5]{ltree2.eps}}
\def\ltc3{\includegraphics[scale=0.8]{ltree3}}
\def\ltd31{\!\!\includegraphics[scale=0.3]{ltree31}}
\def\lte4{\includegraphics[scale=0.5]{ltree4}}

\def\ltde31{\!\!\includegraphics[scale=0.5]{ltree31e.eps}}
\def\ltdn31{\!\!\includegraphics[scale=0.5]{ltree31n.eps}}
\def\ltf41{\includegraphics[scale=0.5]{l41.eps}}
\def\lti4{\includegraphics[scale=0.5]{l411.eps}}
\def\ltg42{\includegraphics[scale=0.5]{l42.eps}}
\def\lth43{\includegraphics[scale=0.5]{l43.eps}}

\def \1e{ \left|1_\mathcal{F}\right>}
\def \2e{\left|\!\!\begin{array}{c}
                    \sfish \\
                   \end{array}\!\!\right>}
\def \3e{\left|\!\!\begin{array}{c}
                    \siceC \\
                 \end{array}\!\!\right>}
\def \s3e{\big|\!\!\begin{array}{c}
                    \xsiceCUp \\
                 \end{array}\!\!\big>}

\def \4e{\left|\!\begin{array}{c}
                    \sfishfish \\
                 \end{array}\!\!\right>}
\nc{\arbreA}{%\kern-0.4ex
\setlength{\unitlength}{.7pt}
\begin{picture}(60,40)(0,0)
\put(30,0){\line(0,1){10}} \put(30,10){\line(-1,1){30}}
\put(30,10){\line(1,1){30}}
\end{picture}}%\kern 0.4ex}

\nc{\arbreAd}{%\kern-0.4ex
\setlength{\unitlength}{.7pt}
\begin{picture}(60,40)(0,0)
\put(30,0){\line(0,1){10}} \put(30,10){\line(-1,1){30}}
\put(30,10){\line(1,1){30}} \put(35,3){\small{x}}
\end{picture}}%\kern 0.4ex}

\nc{\ladderc}{%\kern-0.4ex
\setlength{\unitlength}{.7pt}
\begin{picture}(60,40)(0,0)
\put(30,0){\line(0,1){50}} \put(26,-3){$\bullet$}
\put(26,22){{$\bullet$}} \put(26,45){$\bullet$} \put(38,-3){$a$}
\put(38,22){{$b$}} \put(38,45){$c$}
\end{picture}}%\kern 0.4ex}

\nc{\disp}[1]{\displaystyle{#1}}
\nc{\bin}[2]{ (_{\stackrel{\scs{#1}}{\scs{#2}}})}  %binomial coeff
\nc{\binc}[2]{ \left (\!\! \begin{array}{c} \scs{#1}\\
    \scs{#2} \end{array}\!\! \right )}  %binomial coeff
\nc{\bincc}[2]{  \left ( {\scs{#1} \atop
    \vspace{-.5cm}\scs{#2}} \right )}  %binomial coeff
\nc{\sarray}[2]{\begin{array}{c}#1 \vspace{.1cm}\\ \hline
    \vspace{-.35cm} \\ #2 \end{array}}
\nc{\bs}{\bar{S}} \nc{\ep}{\epsilon}
\nc{\dbigcup}{\stackrel{\bullet}{\bigcup}}
\nc{\la}{\longrightarrow} \nc{\cprod}{\ast} \nc{\rar}{\rightarrow}
\nc{\dar}{\downarrow} \nc{\labeq}[1]{\stackrel{#1}{=}}
\nc{\dap}[1]{\downarrow \rlap{$\scriptstyle{#1}$}}
\nc{\uap}[1]{\uparrow \rlap{$\scriptstyle{#1}$}}
\nc{\defeq}{\stackrel{\rm def}{=}} \nc{\dis}[1]{\displaystyle{#1}}
\nc{\dotcup}{\ \displaystyle{\bigcup^\bullet}\ }
\nc{\sdotcup}{\tiny{ \displaystyle{\bigcup^\bullet}\ }}
\nc{\fe}{\'{e}}
\nc{\hcm}{\ \hat{,}\ } \nc{\hcirc}{\hat{\circ}}
\nc{\hts}{\hat{\shpr}} \nc{\lts}{\stackrel{\leftarrow}{\shpr}}
\nc{\denshpr}{\den{\shpr}}
\nc{\rts}{\stackrel{\rightarrow}{\shpr}} \nc{\lleft}{[}
\nc{\lright}{]} \nc{\uni}[1]{\tilde{#1}} \nc{\free}[1]{\bar{#1}}
\nc{\freea}[1]{\tilde{#1}} \nc{\freev}[1]{\hat{#1}}
\nc{\dt}[1]{\hat{#1}}
\nc{\wor}[1]{\check{#1}}
\nc{\intg}[1]{F_C(#1)}
\nc{\den}[1]{\check{#1}} \nc{\lrpa}{\wr} \nc{\mprod}{\pm}
\nc{\dprod}{\ast_P} \nc{\curlyl}{\left \{ \begin{array}{c} {} \\
{} \end{array}
    \right .  \!\!\!\!\!\!\!}
\nc{\curlyr}{ \!\!\!\!\!\!\!
    \left . \begin{array}{c} {} \\ {} \end{array}
    \right \} }
\nc{\longmid}{\left | \begin{array}{c} {} \\ {} \end{array}
    \right . \!\!\!\!\!\!\!}
\nc{\lin}{\call} \nc{\ot}{\otimes}
\nc{\ora}[1]{\stackrel{#1}{\rar}}
\nc{\ola}[1]{\stackrel{#1}{\la}}%${\Bbb Z}$
\nc{\scs}[1]{\scriptstyle{#1}} \nc{\mrm}[1]{{\rm #1}}
\nc{\margin}[1]{\marginpar{\rm #1}}   %{\rm #1}}
\nc{\dirlim}{\displaystyle{\lim_{\longrightarrow}}\,}
\nc{\invlim}{\displaystyle{\lim_{\longleftarrow}}\,}
\nc{\mvp}{\vspace{0.5cm}}
\nc{\mult}{m}       %multiplication in bialgebra
\nc{\svp}{\vspace{2cm}} \nc{\vp}{\vspace{8cm}}
\nc{\proofbegin}{\noindent{\bf Proof: }}
\nc{\proofend}{$\blacksquare$ \vspace{0.5cm}}
\nc{\sha}{{\mbox{\cyr X}}}  %used to be \cyr
\nc{\ncsha}{{\mbox{\cyr X}^{\mathrm NC}}} \nc{\ncshao}{{\mbox{\cyr
X}^{\mathrm NC,\,0}}}
\nc{\shpr}{\diamond}    %Shuffle product
\nc{\shprc}{\shpr_c}
\nc{\shpro}{\diamond^0}    %Shuffle product
\nc{\shpru}{\check{\diamond}} \nc{\spr}{\cdot}
\nc{\catpr}{\diamond_l} \nc{\rcatpr}{\diamond_r}
\nc{\lapr}{\diamond_a} \nc{\lepr}{\diamond_e} \nc{\sprod}{\bullet}
\nc{\un}{u}                 %unit map in bialgebra
\nc{\vep}{\varepsilon} \nc{\labs}{\mid\!} \nc{\rabs}{\!\mid}
\nc{\hsha}{\widehat{\sha}} \nc{\lsha}{\stackrel{\leftarrow}{\sha}}
\nc{\rsha}{\stackrel{\rightarrow}{\sha}} \nc{\lc}{\lfloor}
\nc{\rc}{\rfloor}
\nc{\sqmon}[1]{\langle #1\rangle}
\nc{\altx}{\Lambda_X}
\nc{\vecT}{\vec{T}}
\nc{\piword}{{\mathfrak P}}

\nc{\mmbox}[1]{\mbox{\ #1\ }}
\nc{\ayb}{\mrm{AYB}} \nc{\mayb}{\mrm{mAYB}} \nc{\cyb}{\mrm{cyb}}
\nc{\ann}{\mrm{ann}} \nc{\Aut}{\mrm{Aut}} \nc{\cabqr}{\mrm{CABQR
}} \nc{\can}{\mrm{can}} \nc{\colim}{\mrm{colim}}
\nc{\Cont}{\mrm{Cont}} \nc{\rchar}{\mrm{char}}
\nc{\cok}{\mrm{coker}} \nc{\dtf}{{R-{\rm tf}}} \nc{\dtor}{{R-{\rm
tor}}}

\nc{\Div}{{\mrm Div}} \nc{\End}{\mrm{End}} \nc{\Ext}{\mrm{Ext}}
\nc{\FG}{\mrm{FG}} \nc{\Fil}{\mrm{Fil}} \nc{\Frob}{\mrm{Frob}}
\nc{\Gal}{\mrm{Gal}} \nc{\GL}{\mrm{GL}} \nc{\Hom}{\mrm{Hom}}
\nc{\hsr}{\mrm{H}} \nc{\hpol}{\mrm{HP}} \nc{\id}{\mrm{id}}
\nc{\im}{\mrm{im}} \nc{\incl}{\mrm{incl}} \nc{\Loday}{\mrm{ABQR}\
} \nc{\length}{\mrm{length}} \nc{\LR}{\mrm{LR}} \nc{\mchar}{\rm
char} \nc{\pmchar}{\partial\mchar} \nc{\map}{\mrm{Map}}
\nc{\MS}{\mrm{MS}} \nc{\OS}{\mrm{OS}} \nc{\NC}{\mrm{NC}}
\nc{\rba}{\rm{Rota-Baxter algebra}\xspace}
\nc{\rbas}{\rm{Rota-Baxter algebras}\xspace}
\nc{\rbw}{\rm{RBW}\xspace}
\nc{\rbws}{\rm{RBWs}\xspace}
\nc{\rbadj}{\rm{RB}\xspace}
\nc{\mpart}{\mrm{part}} \nc{\ql}{{\QQ_\ell}} \nc{\qp}{{\QQ_p}}
\nc{\rank}{\mrm{rank}} \nc{\rcot}{\mrm{cot}} \nc{\rdef}{\mrm{def}}
\nc{\rdiv}{{\rm div}} \nc{\rtf}{{\rm tf}} \nc{\rtor}{{\rm tor}}
\nc{\res}{\mrm{res}} \nc{\SL}{\mrm{SL}} \nc{\Spec}{\mrm{Spec}}
\nc{\tor}{\mrm{tor}} \nc{\Tr}{\mrm{Tr}}
\nc{\mtr}{\mrm{tr}}

\nc{\ab}{\mathbf{Ab}} \nc{\Alg}{\mathbf{Alg}}
\nc{\Bax}{\mathbf{CRB}} \nc{\Algo}{\mathbf{Alg}^0}
\nc{\cRB}{\mathbf{CRB}} \nc{\cRBo}{\mathbf{CRB}^0}
\nc{\RBo}{\mathbf{RB}^0} \nc{\BRB}{\mathbf{RB}}
\nc{\Dend}{\mathbf{DD}} \nc{\bfk}{{\bf k}} \nc{\bfone}{{\bf 1}}
\nc{\base}[1]{{a_{#1}}} \nc{\Cat}{\mathbf{Cat}}

\nc{\DT}{\mathbf{DT}} \nc{\detail}{\marginpar{\bf More detail}
    \noindent{\bf Need more detail!}
    \svp}
\nc{\Diff}{\mathbf{Diff}} \nc{\gap}{\marginpar{\bf
Incomplete}\noindent{\bf Incomplete!!}
    \svp}
\nc{\FMod}{\mathbf{FMod}} \nc{\Int}{\mathbf{Int}}
\nc{\Mon}{\mathbf{Mon}}
\nc{\RB}{\mathbf{RB}} \nc{\remarks}{\noindent{\bf Remarks: }}
\nc{\Rep}{\mathbf{Rep}} \nc{\Rings}{\mathbf{Rings}}
\nc{\Sets}{\mathbf{Sets}} \nc{\bfx}{\mathbf{x}}
\nc{\BA}{{\Bbb A}} \nc{\CC}{{\Bbb C}} \nc{\DD}{{\Bbb D}}
\nc{\EE}{{\Bbb E}} \nc{\FF}{{\Bbb F}} \nc{\GG}{{\Bbb G}}
\nc{\HH}{{\Bbb H}} \nc{\LL}{{\Bbb L}} \nc{\NN}{{\Bbb N}}
\nc{\QQ}{{\Bbb Q}} \nc{\RR}{{\Bbb R}} \nc{\TT}{{\Bbb T}}
\nc{\VV}{{\Bbb V}} \nc{\ZZ}{{\Bbb Z}}

\nc{\cala}{{\mathcal A}} \nc{\calc}{{\mathcal C}}
\nc{\cald}{{\mathcal D}} \nc{\cale}{{\mathcal E}}
\nc{\calf}{{\mathcal F}} \nc{\calg}{{\mathcal G}}
\nc{\calh}{{\mathcal H}} \nc{\cali}{{\mathcal I}}
\nc{\calj}{{\mathcal J}} \nc{\call}{{\mathcal L}}
\nc{\calm}{{\mathcal M}} \nc{\caln}{{\mathcal N}}
\nc{\calo}{{\mathcal O}} \nc{\calp}{{\mathcal P}}
\nc{\calr}{{\mathcal R}} \nc{\calt}{{\mathcal T}}
\nc{\calw}{{\mathcal W}} \nc{\calx}{{\mathcal X}}
\nc{\CA}{\mathcal{A}}

\nc{\frakA}{{\mathfrak A}}
\nc{\fraka}{{\mathfrak a}}
\nc{\frakB}{{\mathfrak B}}
\nc{\frakb}{{\mathfrak b}}
\nc{\frakd}{{\mathfrak d}}
\nc{\frakF}{{\mathfrak F}}
\nc{\frakg}{{\mathfrak g}}
\nc{\frakm}{{\mathfrak m}}
\nc{\frakM}{{\mathfrak M}}
\nc{\frakMo}{{\mathfrak M}^0}
\nc{\frakP}{{\mathfrak P}}
\nc{\frakp}{{\mathfrak p}}
\nc{\frakS}{{\mathfrak S}}
\nc{\frakSo}{{\mathfrak S}^0}
\nc{\fraks}{{\mathfrak s}}
\nc{\os}{\overline{\fraks}}
\nc{\frakT}{{\mathfrak T}}
\nc{\frakTo}{{\mathfrak T}^0}
\nc{\oT}{\overline{T}}
\nc{\frakX}{{\mathfrak X}}
\nc{\frakXo}{{\mathfrak X}^0}
\nc{\frakx}{{\mathbf x}}
\nc{\frakTx}{\frakT}      %All rooted trees, correspond to \ncsha(X)
\nc{\frakTa}{\frakT^a}        % rooted trees for \ncsha(A)
\nc{\frakTxo}{\frakTx^0}   % rooted trees for \ncshao(X)
\nc{\caltao}{\calt^{a,0}}   % rooted trees for \ncshao(A)
\nc{\ox}{\overline{\frakx}}
\nc{\fraky}{{\mathfrak y}}
\nc{\frakz}{{\mathfrak z}}
\nc{\oX}{\overline{X}}

\font\cyr=wncyr10 

\begin{document}

\nc{\redtext}[1]{\textcolor{red}{#1}}

\title[Rota-Baxter algebras and dendriform algebras]
{Rota-Baxter algebras and dendriform algebras}
\author{Kurusch Ebrahimi-Fard}
\address{I.H.\'E.S.
         Le Bois-Marie,
         35, Route de Chartres,
         F-91440 Bures-sur-Yvette,
         France}
         \email{kurusch@ihes.fr}
\author{Li Guo}
\address{Department of Mathematics and Computer Science,
         Rutgers University,
         Newark, NJ 07102}
\email{liguo@newark.rutgers.edu}

%\date{\today}

%\begin{document}

\begin{abstract}
In this paper we study the adjoint functors between the category
of Rota-Baxter algebras and the categories of dendriform dialgebras
and trialgebras.
In analogy to the well-known theory of the adjoint functor between the category of
associative algebras and Lie algebras, we first give an explicit
construction of free Rota-Baxter algebras and then apply it
to obtain universal enveloping Rota-Baxter algebras
of dendriform dialgebras and trialgebras. We further show that free dendriform
dialgebras and trialgebras, as represented by binary planar trees and
planar trees, are canonical subalgebras of free Rota-Baxter algebras.
\end{abstract}

\maketitle

\delete{
\noindent
{\bf Corresponding author}: Li Guo\\
Department of Mathematics and Computer Science,\\
Rutgers University,\\
Newark, NJ 07102, USA\\
liguo@newark.rutgers.edu\\
Phone: (973) 353-5156 ext 30.\\
Fax: (973) 353-5270
}

\noindent
{\bf Keywords: } Rota-Baxter algebra, dendriform algebra,
adjoint functor, universal enveloping algebra, planar tree.

%\tableofcontents

\setcounter{section}{0}

\section{Introduction}
\mlabel{sec:intro}

It is well-known that the natural functor from the category of associative
algebras to that of Lie algebras and the adjoint functor
play a fundamental role in the study of
these algebraic structures and their applications.
This paper establishes a similar relationship between Rota-Baxter algebras
and dendriform dialgebras and dendriform trialgebras by using
free Rota-Baxter algebras.
\medskip

%\subsection{Rota-Baxter algebras}

A Rota-Baxter algebra is an algebra $A$ with a linear endomorphism $R$ satisfying
the {\bf Rota-Baxter equation}:

\begin{equation}
    R(x)R(y) = R\big(R(x)y + xR(y) + \lambda xy\big),\ \forall x,y \in A.
    \mlabel{eq:RB}
\end{equation}
Here $\lambda$ is a fixed element in the base ring and is sometimes denoted by
$-\theta$. This equation was introduced by the mathematician Glen E.
Baxter~\mcite{Ba} in 1960 in his probability study, and was popularized
mainly by the work of Gian-Carlo Rota~\mcite{Ro1, Ro2, Ro3} and his school.

Linear operators satisfying equation
(\mref{eq:RB}) in the context of Lie algebras were introduced
independently by Belavin and Drinfeld \mcite{B-D}, and
Semenov-Tian-Shansky~\mcite{STS1} in the 1980s and were related to
solutions, called $r$-matrices, of the (modified) classical
Yang-Baxter equation, named after the physicists Chen-ning Yang
and Rodney Baxter.
Recently, there have been several interesting developments of Rota-Baxter
algebras in theoretical physics and
mathematics, including quantum field
theory~\mcite{C-K1,C-K2}, Yang-Baxter
equations~\mcite{Ag1,Ag2,Ag3}, shuffle
products~\mcite{shuf,G-K1,G-K2},
operads~\mcite{A-L,EF1,prod,Le1,Le2,Le3}, Hopf
algebras~\mcite{A-G-K-O,shuf,EMP07}, combinatorics~\mcite{Gu2} and
number theory~\mcite{shuf,mzv,Gu5,G-Z,MP1,MP2,zhao}. The most prominent of
these is the work~\mcite{C-K1,C-K2} of Connes and Kreimer in their Hopf algebraic
approach to renormalization theory in perturbative quantum field
theory, continued in a series
of papers
\mcite{E-G-G-V,mat,EGfields06,E-G-K2,E-G-K3,egm2006,ek2005,
em2006,EMP07}.
\smallskip

%\subsection{Dendriform dialgebras and trialgebras}

A dendriform dialgebra is a module $D$ with two binary operations
$\prec$ and $\succ$ that satisfy three relations between them (see Eq.~(\mref{eq:dia})).
This concept was introduced by Loday~\mcite{Lo1} in 1995
with motivation from algebraic $K$-theory, and was further
studied in connection with several areas in mathematics and
physics, including operads~\mcite{Lo2}, homology~\mcite{Fra1,Fra2},
Hopf algebras~\mcite{Ch,Hol2,L-R2,Ron,KDF2007}, Lie and Leibniz algebras~\mcite{Fra2},
combinatorics~\mcite{A-S1,A-S2,Fo,L-R1}, arithmetic~\mcite{Lo3} and
quantum field theory~\mcite{Fo,Hol1}.

A few years later Loday and Ronco defined dendriform trialgebras
in their study~\mcite{L-R2} of polytopes and Koszul duality.
Such a structure is a module $T$ equipped with binary operations
$\prec,\succ$ and
$\spr$ that satisfy seven relations that will be recalled in Eq.~(\mref{eq:tri}).

The dendriform dialgebra and trialgebra share the property that the sum of
the binary operations $\prec+\succ$ (for dialgebra) or $\prec+\succ+\, \spr$
(for trialgebra) is associative. Other dendriform algebra structures
have the similar property of ``splitting associativity" in the sense that
an associative
product decomposes into a linear combination of several binary operations.
Many such structures have been obtained lately, such as
the quadri-algebra of Loday and Aguiar~\mcite{A-L} and
the ennea- and NS-algebra of Leroux~\mcite{Le1,Le2}. In \mcite{prod}
(see also \mcite{Lo4}),
we showed how these more complex structures, equipped with large
numbers of compositions and relations, can be derived from an
operadic point of view in terms of products.
Further examples and developments can be found in~\mcite{unit,Lo2}.
\smallskip
%%%
%%%%%%%%%%%%%%%%%%%%%%%%%%%%%%%%%%%%%%%%%%%%%%%%%%%%%%%%%%%
%%%

%\subsection{The connection}
The first link between Rota-Baxter algebras and dendriform
algebras was given by Aguiar~\mcite{Ag1} who showed that a
Rota-Baxter algebra of weight $\lambda=0$ carries a dendriform
dialgebra structure, resembling the Lie algebra structure on an associative
algebra. This has been extended to
further connections between linear operators and dendriform type
algebras~\mcite{EF1,Le2,A-L,prod}, in particular to dendriform trialgebras by
the first named author. See Theorem~\mref{thm:EFs} for details.

Consequently, there are natural functors from the category of
Rota-Baxter algebras of weight $\lambda$ to the categories of
dendriform dialgebras and trialgebras. We study the adjoint functors in this paper.
\smallskip

%\subsection{Outline of the paper}

As a preparation, we first construct in Section~\mref{sec:nonua}
free Rota-Baxter algebras (Theorem~\mref{thm:freeao}) which play a central role in the study
of the adjoint functors. This is in analogy to the central role played by
the free associative algebras in the study of the adjoint functor
from the category of Lie algebras to the category of associative
algebras.
As we will see, free Rota-Baxter algebras can be defined in various
generalities, such as over a set or over another algebra, in various
contexts, such as unitary or nonunitary algebras, and they can be
constructed in various terms, such as by words or by trees, either explicitly
or recursively. For the purpose of our application to adjoint functors,
we only consider a special case of free Rota-Baxter algebras, namely
free nonunitary Rota-Baxter algebras $\ncshao(A)$ generated by another algebra $A$ that
possesses a basis over the base ring. Further studies of free Rota-Baxter algebras can be found in~\mcite{A-M,free,EMP07,Gu6,GK3,G-S}.

Then in Section \mref{sec:adj},
we use these free Rota-Baxter algebras to
obtain adjoint functors of the functors from Rota-Baxter algebras
to dendriform dialgebras (Theorem~\mref{thm:envdend}) and trialgebras (Theorem~\mref{thm:env}) by proving the existence of the corresponding universal enveloping Rota-Baxter algebras. In the case of dendriform trialgebras,
let $D=(D,\prec,\succ,\spr)$ be a dendriform trialgebra. Let $\ncshao(D)$ be the free nonunitary
Rota-Baxter algebra over the nonunitary algebra $(D,\spr)$ constructed in
Theorem~\mref{thm:freeao}. Let $I_R$ be
a suitable Rota-Baxter ideal of $\ncshao(D)$ generated by relations from $\prec$ and $\succ$. Theorem~\mref{thm:envdend} shows that
the quotient Rota-Baxter algebra $\ncshao(D)/I_R$ is the universal enveloping Rota-Baxter algebra of $D$ in the sense of Definition~\mref{de:env}.

The special case of free dendriform algebras is considered in
Section~\mref{sec:dfree} where we realize the free dendriform
dialgebra and trialgebra of Loday and Loday-Ronco
in terms of decorated planar rooted trees as canonical subalgebras of free Rota-Baxter algebras.

\medskip

\noindent
{\bf Notations:}
In this paper, $\bfk$ is a commutative unitary ring which will be
further assumed to be a field in Sections~\mref{sec:adj} and
\mref{sec:dfree}. Let $\Alg$ be the category of unitary
$\bfk$-algebras $A$ whose unit is identified with the unit
$\bfone$ of $\bfk$ by the structure homomorphism $\bfk\to A$. Let
$\Algo$ be the category of nonunitary $\bfk$-algebras. Similarly
let $\RB_\lambda$ (resp. $\RBo_\lambda$) be the category of
unitary (resp. nonunitary) Rota-Baxter $\bfk$-algebras of weight
$\lambda$. The subscript $\lambda$ will be suppressed if there is
no danger of confusion.
\medskip

%%%
%%%%%%%%%%%%%%%%%%%%%%%%%%%%%%%%%%%%%%%%%%%%%%%%%%%%%%%%%%%
%%%

\noindent
{\bf Acknowledgements:}
We thank M. Aguiar, J.-L. Loday and M. Ronco for helpful
discussions. The first named author was supported by a Ph.D. grant
from the Ev. Studienwerk e.V., and would like to thank the people
at the Theory Department of the Physics Institute at Bonn
University for encouragement and help. The second named author
acknowledges support from NSF grant DMS 0505643 and
a Research Council grant from the
Rutgers University. Both authors acknowledge the warm hospitality
of I.H.\'E.S. (LG) and L.P.T.H.E. (KEF) where this work was
completed.

%%%
%%%%%%%%%%%%%%%%%%%%%%%%%%%%%%%%%%%%%%%%%%%%%%%%%%%%%%%%%%%
%%%

\section{Free nonunitary Rota-Baxter algebras on an algebra}
\mlabel{sec:nonua}

We now construct free nonunitary \rbas over another nonunitary
algebra. Other than its theoretical significance, our main purpose
is for the application in later sections to
study universal enveloping \rbas of dendriform dialgebras and
trialgebras.
The reader can regard such free \rbas over another algebra as
the Rota-Baxter analog of the tensor algebra over a module.
It is well-known that such tensor algebras are essential in the
study of enveloping algebras of Lie algebras~\mcite{Reu}. Because of the
nonunitariness of Lie algebras, it is the free nonunitary, instead of
unitary, associative algebras
that are used in the study of the adjoint functor from Lie algebra to
associative algebras. For the similar reason, free nonunitary Rota-Baxter
algebras are convenient in the study of the adjoint functor from dendriform
algebras to Rota-Baxter algebras.
As remarked earlier, other cases of free Rota-Baxter algebras are considered
elsewhere~\mcite{free}.

Let $B$ be a nonunitary $\bfk$-algebra. Recall~\mcite{G-K1,G-K2} that
a free nonunitary \rba over $B$ is a nonunitary \rba $\ncshao(B)$ with a Rota-Baxter
operator $R_B$ and a nonunitary algebra homomorphism $j_B: B\to \ncshao(B)$ such that,
for any nonunitary \rba $A$ and any nonunitary algebra
homomorphism $f:B\to A$, there is a unique nonunitary \rba homomorphism
$\free{f}: \ncshao(B)\to A$ such that $\free{f}\circ j_B=f$.
$$ \xymatrix{ B \ar[rr]^{j_B}\ar[drr]^{f} && \ncshao(B) \ar[d]_{\free{f}} \\
&& A}
$$

We assume that the nonunitary algebra $B$ possesses a basis over the
base ring $\bfk$. This is no restriction if the base ring is a field
as is customarily taken to be the case in the study of dendriform
algebras/operads and therefore in our later sections.

We first display a $\bfk$-basis of the free \rba in terms of words
in \S~\mref{ss:base}. The product on the free
\rba is given in~\mref{ss:prodao} and the universal property of
the free \rba is proved in~\mref{ss:proof}.

\subsection{A basis of a free Rota-Baxter algebra as words}
\mlabel{ss:base}
Let $B$ be a nonunitary $\bfk$-algebra with a $\bfk$-basis $X$.
We first display a $\bfk$-basis $\frakX_\infty$ of $\ncshao(B)$ in terms of
words from the alphabet set $X$.

Let $\lc$ and $\rc$ be symbols, called brackets,
and let $X'=X\cup \{\lc,\rc\}$.
Let $M(X')$ be the free semigroup generated by $X'$.

\begin{defn}
Let $Y,Z$ be two subsets of $M(X')$. Define the {\bf alternating product}
of $Y$ and $Z$ to be
\allowdisplaybreaks{
\begin{eqnarray}
\altx(Y,Z)&=&\Big( \bigcup_{r\geq 1} \big (Y\lc Z\rc \big)^r \Big) \bigcup
    \Big(\bigcup_{r\geq 0} \big (Y\lc Z\rc \big)^r  Y\Big) \notag \\
&& \bigcup \Big( \bigcup_{r\geq 1} \big( \lc Z\rc Y \big )^r \Big)
 \bigcup \Big( \bigcup_{r\geq 0} \big (\lc Z\rc Y\big )^r \lc Z\rc \Big).
\mlabel{eq:wordsao}
\end{eqnarray}}
\mlabel{de:alt}
\end{defn}

We construct a sequence $\frakX_n$ of subsets of $M(X')$ by the following
recursion. Let $\frakX_0=X$ and, for $n\geq 0$, define
\allowdisplaybreaks{
\begin{equation}
\frakX_{n+1}=\altx(X,\frakX_n). \notag
\end{equation}
More precisely,
\begin{eqnarray}
\frakX_{n+1}&=& \Big( \bigcup_{r\geq 1} \big (X\lc \frakX_{n}\rc\big )^r \Big) \bigcup
    \Big(\bigcup_{r\geq 0} \big (X\lc \frakX_{n}\rc\big )^r  X\Big)
    \notag \\
&& \bigcup \Big( \bigcup_{r\geq 1} \big (\lc \frakX_{n}\rc X\big )^r \Big)
    \bigcup \Big( \bigcup_{r\geq 0} \big( \lc \frakX_{n}\rc X \big)^r
    \lc \frakX_{n-1}\rc \Big). \mlabel{eq:x1ao}
\end{eqnarray}
Further, define
\begin{eqnarray}
\frakX_\infty &=& \bigcup_{n\geq 0} \frakX_n = \dirlim \frakX_n. \mlabel{eq:x3ao}
\end{eqnarray}}
Here the second equation in Eq. (\mref{eq:x3ao}) follows since
$\frakX_1\supseteq \frakX_0$ and, assuming $\frakX_n\supseteq \frakX_{n-1}$,
we have
$$\frakX_{n+1}=\altx(X,\frakX_n) \supseteq \altx(X,\frakX_{n-1})
    \supseteq \frakX_n.$$
\begin{defn}
A word in $\frakX_\infty$ is called a {\bf (strict) Rota-Baxter (bracketed)
word (RBWs)}.
\mlabel{de:rbw}
\end{defn}

A similar concept of parenthesized words has appeared in the work of
Kreimer~\mcite{Kr1} to represent Hopf algebra structure on Feynman
diagrams in pQFT, with a different set of restrictions on the
words. We use the brackets $\lc$ and $\rc$ instead of $($ and $)$ to avoid
confusion with the
usual meaning of parentheses.

The verification of the following properties of RBWs are quite easy and is left to the reader.
\begin{lemma}
\begin{enumerate}
\item
For each $n\geq 1$, the union of
$\frakX_n=\altx(X,\frakX_{n-1})$ expressed in Eq.(\mref{eq:x1ao}) is
disjoint:
\allowdisplaybreaks{
\begin{eqnarray}
\frakX_n & =&
    \Big( \dbigcup_{r\geq 1} \big (X\lc \frakX_{n-1}\rc\big )^r \Big) \dbigcup
    \Big(\dbigcup_{r\geq 0} \big (X\lc \frakX_{n-1}\rc\big )^r  X\Big) \notag\\
 && \dbigcup \Big( \dbigcup_{r\geq 1} \big (\lc \frakX_{n-1}\rc X\big )^r \Big)
    \dbigcup \Big( \dbigcup_{r\geq 0} \big( \lc \frakX_{n-1}\rc X \big)^r
    \lc \frakX_{n-1}\rc \Big).
\mlabel{eq:words2}
\end{eqnarray}}
\mlabel{it:disjoint}
\item
We further have the disjoint union
\allowdisplaybreaks{
\begin{eqnarray}
\frakX_\infty & =&
    \Big( \dbigcup_{r\geq 1} \big (X\lc \frakX_{\infty}\rc\big )^r \Big) \dbigcup
    \Big(\dbigcup_{r\geq 0} \big (X\lc \frakX_{\infty}\rc\big )^r  X\Big) \notag\\
 && \dbigcup \Big( \dbigcup_{r\geq 1} \big (\lc \frakX_{\infty}\rc X\big )^r \Big)
    \dbigcup \Big( \dbigcup_{r\geq 0} \big( \lc \frakX_{\infty}\rc X \big)^r
    \lc \frakX_{\infty}\rc \Big).
\mlabel{eq:words3}
\end{eqnarray}}
\mlabel{it:disjointi}
\item
Every RBW $\frakx\neq \bfone$ has a unique decomposition
\begin{equation}
 \frakx=\frakx_1 \cdots  \frakx_b,
\mlabel{eq:st}
\end{equation}
where $\frakx_i$, $1\leq i\leq b$, is alternatively in $X$ or in $\lc \frakX_\infty\rc$.
This decomposition will be called the {\bf standard decomposition}
of $\frakx$.
\mlabel{it:st}
\end{enumerate}
\mlabel{lem:ex}
\end{lemma}
For a \rbw $\frakx$ in ${\frakX}_\infty$ with standard decomposition
$\frakx_1 \cdots  \frakx_b$,
we define $b$ to be the {\bf breadth} $b(\frakx)$ of $\frakx$, we define
the {\bf head}
$h(\frakx)$ of $\frakx$ to be 0 (resp. 1) if $\frakx_1$ is in $X$
(resp. in $\lc \frakX_\infty \rc$). Similarly define the {\bf tail}
$t(\frakx)$ of $\frakx$ to be 0 (resp. 1) if $\frakx_b$ is in $X$
(resp. in $\lc \frakX_\infty \rc$).
In terms of the decomposition~(\mref{eq:words2}),
the head, tail and breadth of a word $\frakx$
are given in the following table.
\begin{center}
\begin{tabular}{c|c|c|c|c}
$\frakx$ & $(X\lc \frakX_{n-1}\rc)^{r}$
&$(X \lc \frakX_{n-1}\rc)^{ r}  X$ &
$(\lc \frakX_{n-1}\rc  X)^{ r} $ &
$(\lc \frakX_{n-1}\rc  X)^{ r}  \lc \frakX_{n-1}\rc$ \\ \hline
$h(\frakx)$& $0$ & 0& 1& 1 \\
$t(\frakx)$ & $1$ & 0 & 0 & 1\\
$b(\frakx)$& $2r$ &$2r+1$& $2r$& $2r+1$
\end{tabular}
\end{center}
Finally, define the {\bf depth} $d(\frakx)$ to be
$$ d(\frakx)=\min \{n\ \big |\ \frakx\in \frakX_n \}.$$
So, in particular, the depth of elements in $X$ is 0 and depth of elements
in $\lc X\rc$ is one.
\begin{exam} For
$x_1,x_2, x_3\in X$, the word $\lc\lc x_1\rc x_2\rc x_3$ has
head 1, tail 0, breadth 2 and depth 2.
\end{exam}

\subsection{The product in a free Rota-Baxter algebra}
\mlabel{ss:prodao}
Let
$$\ncshao(B)=\bigoplus_{\frakx\in \frakX_\infty} \bfk \frakx.$$
We now define a product $\shpr$ on $\ncshao(B)$ by defining
$\frakx\shpr \frakx'\in \ncshao(B)$ for $\frakx,\frakx'\in \frakX_\infty$ and then
extending bilinearly.
Roughly speaking, the product of $\frakx$ and $\frakx'$ is defined
to be the concatenation whenever $t(\frakx)\neq h(\frakx')$. When
$t(\frakx)=h(\frakx')$, the product is defined by the product in
$B$ or by the Rota-Baxter relation in Eq.~(\mref{eq:shprod0}).

To be precise, we use induction on the sum $n:=d(\frakx)+d(\frakx')$
of the depths of $\frakx$ and $\frakx'$.
Then $n\geq 0$.
If $n=0$, then $\frakx,\frakx'$ are in $X$ and so are in $B$ and we
define $\frakx\shpr \frakx'=\frakx \spr \frakx'\in B \subseteq \ncshao(B)$.
Here $\spr$ is the product in $B$.

Suppose $\frakx\shpr \frakx'$ have been defined for all $\frakx,\frakx'\in
\frakX_\infty$ with $n\geq k\geq 0$ and let
$\frakx, \frakx'\in \frakX_\infty$ with $n=k+1$.

First assume the breadth $b(\frakx)=b(\frakx')=1$. Then
$\frakx$ and $\frakx'$ are in $X$ or $\lc \frakX_\infty\rc$. We accordingly
define
\begin{equation}
\frakx\shpr \frakx'=\left \{ \begin{array}{ll}
\frakx \spr \frakx', & {\rm if\ } \frakx,\frakx'\in X,\\
\frakx \frakx', & {\rm if\ } \frakx\in X, \frakx'\in \lc \frakX_\infty\rc,\\
\frakx \frakx', & {\rm if\ } \frakx\in \lc \frakX_\infty\rc, \frakx'\in X,\\
\lc \lc \ox\rc \shpr \ox'\rc +\lc \ox \shpr \lc \ox'\rc \rc
+\lambda \lc \ox \shpr \ox' \rc, & {\rm if\ } \frakx=\lc \ox\rc,
\frakx'=\lc \ox'\rc \in \lc \frakX_\infty \rc.
\end{array} \right .
\mlabel{eq:shprod0}
\end{equation}
Here the product in the first case is the product in $B$, in the second and
third case are by concatenation and in the fourth case is by the induction
hypothesis since for the three products on the right hand side we have
\begin{eqnarray*}
d(\lc\ox \rc)+ d(\ox')
&=& d(\lc \ox \rc)+d(\lc \ox' \rc)-1
= d(\frakx)+d(\frakx')-1,\\
d(\ox)+d(\lc \ox'\rc) &=& d(\lc \ox \rc)+d(\lc \ox'\rc)-1
= d(\frakx)+ d(\frakx')-1,\\
d(\ox)+ d(\ox') &=& d(\lc \ox \rc)-1+ d(\lc \ox' \rc)-1
= d(\frakx)+d(\frakx')-2
\end{eqnarray*}
which are all less than or equal to $k$.

Now assume $b(\frakx)>1$ or $b(\frakx')>1$. Let
$\frakx=\frakx_1\cdots\frakx_b$ and $\frakx'=\frakx'_1\cdots\frakx'_{b'}$
be the standard decompositions from Lemma~\mref{lem:ex}. We then define
\begin{equation}
\frakx \shpr \frakx'= \frakx_1\cdots \frakx_{b-1}(\frakx_b\shpr \frakx'_1)\,
    \frakx'_{2}\cdots \frakx'_{b'}
\end{equation}
where $\frakx_b\shpr \frakx'_1$ is defined by Eq.~(\mref{eq:shprod0}) and
the rest is given by concatenation. The concatenation is well-defined since by
Eq.~(\mref{eq:shprod0}), we have $h(\frakx_b)=h(\frakx_b\shpr \frakx'_1)$
and $t(\frakx'_1)=t(\frakx_b\shpr \frakx'_1)$. Therefore,
$t(\frakx_{b-1})\neq h(\frakx_b\shpr \frakx'_1)$ and
$h(\frakx'_2)\neq t(\frakx_b\shpr \frakx'_1)$.

\medskip

We record the following simple properties of $\shpr$ for later applications.
\begin{lemma} Let $\frakx,\frakx'\in \frakX_\infty$. We have the following
statements.
\begin{enumerate}
\item $h(\frakx)=h(\frakx\shpr \frakx')$ and $t(\frakx')=t(\frakx\shpr \frakx')$.
\mlabel{it:mat0}
\item If $t(\frakx)\neq h(\frakx')$, then
$\frakx \shpr \frakx' =\frakx \frakx'$ (concatenation).
\mlabel{it:mat1}
\item If $t(\frakx)\neq h(\frakx')$, then for any $\frakx''\in \frakX_\infty$,
$$(\frakx\frakx')\shpr \frakx'' =\frakx(\frakx' \shpr \frakx''), \quad
\frakx''\shpr (\frakx \frakx') =(\frakx'' \shpr \frakx) \frakx'.$$
\mlabel{it:mat2}
\end{enumerate}
\mlabel{lem:match}
\end{lemma}

Extending $\shpr$ bilinearly, we obtain
a binary operation
$$ \ncshao (B)\otimes \ncshao(B) \to \ncshao(B).$$
For $\frakx\in \frakX_\infty$, define
\begin{equation}
R_B(\frakx)=\lc \frakx \rc.
\mlabel{eq:RBop}
\end{equation}
Obviously $\lc \frakx \rc$ is again in $\frakX_\infty$. Thus $R_B$ extends to
a linear operator $R_B$ on $\ncshao(B)$.
Let
$$j_X:X\to \frakX_\infty \to \ncshao(B)$$
be the natural injection which extends to an algebra injection
\begin{equation}
j_B: B \to \ncshao(B).
\mlabel{eq:jo}
\end{equation}

The following is our first main result which will be proved in the next subsection.
\begin{theorem}
Let $B$ be a nonunitary $\bfk$-algebra with a $\bfk$-basis $X$.
\begin{enumerate}
\item
The pair $(\ncshao(B),\shpr)$ is a nonunitary associative algebra.
\mlabel{it:alg}
\item
The triple $(\ncshao(B),\shpr,R_B)$ is a nonunitary \rba of weight $\lambda$.
\mlabel{it:RB}
\item
The quadruple $(\ncshao(B),\shpr,R_B,j_B)$ is the free nonunitary \rba of
weight $\lambda$ on
the algebra $B$.
\mlabel{it:free}
\end{enumerate}
\mlabel{thm:freeao}
\end{theorem}

The following corollary of the theorem will be used later in the paper.
\begin{coro}
Let $V$ be a $\bfk$-module and let $T(V)=\bigoplus_{n\geq 1} V^{\ot n}$
be the tensor algebra over $V$. Then $\ncshao(T(V))$, together with the
natural injection $i_V: V\to T(V) \xrightarrow{j_{T(V)}} \ncshao(T(V))$,
is a free nonunitary Rota-Baxter algebra over $V$, in the sense that,
for any nonunitary Rota-Baxter algebra $A$ and $\bfk$-module map
$f: V\to A$ there is a unique nonunitary Rota-Baxter algebra homomorphism
$\freev{f}: \ncshao(T(V)) \to A$ such that $k_V \circ \free{f} = f$.
\mlabel{co:vecfree}
\end{coro}
\begin{proof}
The maps in the corollary and in this proof are organized in the following
diagram
$$ \xymatrix{ T \ar[rr]^{k_V} \ar[dd]_{f} \ar[rrdd]^(.4){i_V}
&& T(V) \ar[dd]^{j_{T(V)}} \ar[lldd]_(.7){\free{f}} \\
&& \\
A && \ncshao(T(V)) \ar[ll]_{\freev{f}}}$$

For the given $\bfk$-module $V$, note that $T(V)$, together with the natural
injection $k_V: V\to T(V)$, is the free nonunitary
$\bfk$-algebra over $V$. So for the given $\bfk$-algebra $A$ and $\bfk$-module
map $f: V\to A$, there is a unique nonunitary $\bfk$-algebra homomorphism
$\freea{f}: T(V) \to A$ such that $\freea{f} \circ k_V=f$. Then by the
universal property of the free Rota-Baxter algebra $\ncshao(T(V))$, there is
a unique $\free{\freea{f}}: \ncshao(T(V)) \to A$ such that
$ \free{\freea{f}}\circ j_{T(V)}=\freea{f}$. Since $i_V=j_{T(V)}\circ k_V$,
we have
$ \free{\freea{f}} i_V = \freea{f} \circ k_V=f$.
So we have proved the existence of $\freev{f}=\free{\freea{f}}.$

For the uniqueness of $\freev{f}$. Suppose there is another
$\freev{f}':\ncshao(T(V)) \to A$ such that $ \freev{f}'\circ i_V =f$.
Then we have
$$ \freev{f}' \circ j_{T(V)} \circ k_V = \freev{f}'\circ i_V=f
    = \freev{f} \circ i_V = \freev{f} \circ j_{T(V)} \circ k_V.$$
By the universal property of the free algebra $T(V)$, we have
$\freev{f}'\circ j_{T(V)} = \freev{f}\circ j_{T(V)}$.
Then by the universal property of the free Rota-Baxter algebra
$\ncshao(T(V))$, we have
$\freev{f}'=\freev{f}$, as needed.
\end{proof}

\subsection{The proof of Theorem~\mref{thm:freeao}}
\mlabel{ss:proof}
\begin{proof}
\mref{it:alg}. We just need to verify the associativity. For this we only need to verify
\begin{equation}
 (\frakx'\shpr \frakx'')\shpr \frakx''' =\frakx'\shpr(\frakx'' \shpr \frakx''')
\mlabel{eq:assx}
\end{equation}
for $\frakx',\frakx'',\frakx'''\in \frakX_\infty$.
We will do this by induction on the sum of the depths
$n:=d(\frakx')+d(\frakx'')+d(\frakx''')$. If $n=0$, then
all of $\frakx',\frakx'',\frakx'''$ have depth zero and so are
in $X$. In this case the product $\shpr$ is given by the product $\spr$
in $B$ and so is associative.

Assume the associativity holds for $n\leq k$ and assume that
$\frakx',\frakx'',\frakx'''\in \frakX_\infty$ have
$n=d(\frakx')+d(\frakx'')+d(\frakx''')=k+1.$

If $t(\frakx')\neq h(\frakx'')$, then by Lemma~\mref{lem:match},
$$ (\frakx' \shpr \frakx'') \shpr \frakx'''=(\frakx'\frakx'')\shpr \frakx'''
= \frakx' (\frakx'' \shpr \frakx''') =\frakx'\shpr (\frakx''\shpr \frakx''').$$
Similarly if $t(\frakx'')\neq h(\frakx''')$.

Thus we only need to verify the associativity when
$t(\frakx')=h(\frakx'')$ and $t(\frakx'')=h(\frakx''')$.
We next reduce the breadths of the words.

\begin{lemma}
If the associativity
$$(\frakx' \shpr \frakx'')\shpr \frakx'''=
 \frakx'\shpr (\frakx'' \shpr \frakx''') $$
holds for all $\frakx', \frakx''$ and $\frakx'''$ in $\frakX_\infty$ of breadth one, then
it holds for all $\frakx', \frakx''$ and $\frakx'''$ in $\frakX_\infty$.
\mlabel{lem:ell}
\end{lemma}

\begin{proof}
We use induction on the sum of breadths
$m:=b(\frakx')+b(\frakx'')+b(\frakx''')$.
Then $m\geq 3$. The case when $m=3$ is the assumption of the lemma.
Assume the associativity holds for
$3\leq m \leq j$ and take $\frakx',
\frakx'',\frakx'''\in \frakX_\infty$ with
$m = j+1.$
Then $j+1\geq 4$. So at least one of $\frakx',\frakx'',\frakx'''$ have
breadth greater than or equal to 2.

First assume $b(\frakx')\geq 2$. Then $\frakx'=\frakx'_1\frakx'_2$
with $\frakx'_1,\, \frakx'_2\in \frakX_\infty$ and
$t(\frakx'_1)\neq h(\frakx'_2)$.
Thus
\allowdisplaybreaks{\begin{eqnarray*}
 (\frakx'\shpr \frakx'') \shpr \frakx'''&=&
((\frakx'_1\frakx'_2)\shpr \frakx'')\shpr \frakx'''\\
&=& (\frakx'_1 (\frakx'_2 \shpr \frakx''))\shpr \frakx'''
 \quad {\rm by\ Lemma~\mref{lem:match}.\mref{it:mat2}}\\
& =& \frakx'_1 ((\frakx'_2 \shpr \frakx'') \shpr \frakx''')
\quad {\rm by\ Lemma~\mref{lem:match}.\mref{it:mat0}\ and\ \mref{it:mat2}}.
\end{eqnarray*}}
Similarly, \allowdisplaybreaks{
\begin{eqnarray*}
 \frakx'\shpr (\frakx'' \shpr \frakx''')&=&
(\frakx'_1\frakx'_2)\shpr (\frakx''\shpr \frakx''')\\
&=& \frakx'_1 (\frakx'_2 \shpr (\frakx''\shpr \frakx''')).
\end{eqnarray*}}
Thus $$ (\frakx'\shpr \frakx'') \shpr \frakx'''=
  \frakx'\shpr (\frakx'' \shpr \frakx''')$$
whenever
$$ (\frakx'_2 \shpr \frakx'') \shpr \frakx'''=
\frakx'_2 \shpr (\frakx''\shpr \frakx''')$$
which follows from the induction hypothesis.

A similar proof works if $b(\frakx''')\geq 2.$

Finally if $b(\frakx'')\geq 2$, then $\frakx''=\frakx''_1\frakx''_2$
with $\frakx''_1,\,\frakx''_2\in \frakX_\infty$ and $t(\frakx''_1)\neq
h(\frakx''_2)$. So using Lemma~\mref{lem:match} repeatedly, we
have \allowdisplaybreaks{
\begin{eqnarray*}
(\frakx' \shpr \frakx'')\shpr \frakx'''&=&
(\frakx' \shpr (\frakx''_1 \frakx''_2)) \shpr \frakx''' \\
&=& ((\frakx' \shpr \frakx''_1)\frakx''_2)\shpr \frakx'''
    \quad {\rm by\ Lemma~\mref{lem:match}.\mref{it:mat0}\ and\ \mref{it:mat2}}\\
&=& (\frakx'\shpr \frakx''_1)(\frakx''_2 \shpr \frakx''')
    \quad {\rm by\ Lemma~\mref{lem:match}.\mref{it:mat0}\ and\ \mref{it:mat2}}
\end{eqnarray*}}
In the same way, we have
$$(\frakx'\shpr \frakx''_1)(\frakx''_2 \shpr \frakx''')
= \frakx'\shpr (\frakx'' \shpr \frakx''').$$
This again proves the associativity.
\end{proof}

To summarize, our proof of the associativity
has been reduced to the special case when
$\frakx',\frakx'',\frakx'''\in \frakX_\infty$ are chosen so that
\begin{enumerate}
\item
$n:= d(\frakx')+d(\frakx'')+d(\frakx''')=k+1\geq 1$ with the assumption that
the associativity holds when $n\leq k$.
\mlabel{it:sp1}
\item
the elements are of breadth one and
\mlabel{it:sp2}
\item
$t(\frakx')=h(\frakx'')$ and $t(\frakx'')=h(\frakx''')$.
\mlabel{it:sp3}
\end{enumerate}
By item \mref{it:sp2}, the head and tail of each of the elements are the same.
Therefore by item \mref{it:sp3}, either all the three elements are in $X$
or they are all in $\lc \frakX_\infty \rc$.
If all of $\frakx',\frakx'',\frakx'''$ are in $X$,
then as already shown, the associativity follows from the associativity in $B$.

So it remains to consider $\frakx',\frakx'',\frakx'''$ all in $\lc
\frakX_\infty \rc$.
Then $\frakx'=\lc
\ox'\rc, \frakx''=\lc \ox'' \rc, \frakx'''=\lc \ox'''\rc$ with
$\ox',\ox'',\ox'''\in \frakX_\infty$. Using Eq.~(\mref{eq:shprod0}) and bilinearity
of the product $\shpr$,
we have
\allowdisplaybreaks{\begin{eqnarray*} (\frakx'\shpr \frakx'')\shpr
\frakx''&=& \big \lc \lc \ox'\rc \shpr \ox '' +\ox'\shpr
\lc\ox''\rc
    +\lambda\ox'\shpr \ox'' \big \rc \shpr \lc \ox'''\rc \\
&=& \lc\lc \ox'\rc \shpr \ox''\rc \shpr \lc\ox'''\rc
    + \lc\ox'\shpr \lc \ox''\rc \rc\shpr \lc \ox'''\rc
    +\lambda \lc \ox'\shpr \ox''\rc \shpr \lc\ox'''\rc \\
&=&  \lc\lc\lc \ox'\rc\shpr \ox''\rc\shpr \ox''' \rc
    + \lc\big(\lc\ox'\rc \shpr \ox''\big) \shpr \lc\ox'''\rc\rc
    +\lambda \lc\big(\lc\ox'\rc \shpr\ox''\big)\shpr \ox'''\rc\\
&& + \lc\lc\ox'\shpr\lc\ox''\rc\rc \shpr \ox'''\rc
    + \lc\big(\ox'\shpr\lc \ox''\rc\big) \shpr\lc \ox'''\rc\rc
    +\lambda \lc\big(\ox'\shpr \lc \ox''\rc \big) \shpr \ox'''\rc \\
&& + \lambda \lc \lc \ox'\shpr \ox''\rc\shpr \ox'''\rc
    +\lambda \lc \big(\ox'\shpr \ox''\big)\shpr \lc \ox'''\rc \rc
    + \lambda^2 \lc \big(\ox'\shpr \ox''\big) \shpr \ox'''\rc.
\end{eqnarray*}}
Applying the induction hypothesis in $n$ to the fifth term $\big
(\ox'\shpr\lc \ox''\rc\big) \shpr\lc \ox'''\rc$ and then use
Eq.~(\mref{eq:shprod0}) again, we have \allowdisplaybreaks{
\begin{eqnarray*}
(\frakx'\shpr \frakx'')\shpr \frakx''
&=& \lc\lc\lc \ox'\rc\shpr \ox''\rc\shpr \ox''' \rc
    + \lc\big(\lc\ox'\rc \shpr \ox''\big) \shpr \lc\ox'''\rc\rc
    +\lambda \lc\big(\lc\ox'\rc \shpr\ox''\big)\shpr \ox'''\rc\\
&& + \lc\lc\ox'\shpr\lc\ox''\rc\rc \shpr \ox'''\rc
    + \lc\ox'\shpr\lc\lc\ox''\rc\shpr\ox'''\rc\rc
    + \lc\ox' \shpr \lc\ox''\shpr \lc \ox'''\rc\rc\rc \\
&&  +\lambda \lc \ox' \shpr \lc\ox''\shpr \ox'''\rc\rc
    +\lambda \lc\big(\ox'\shpr \lc \ox''\rc \big) \shpr \ox'''\rc \\
&& + \lambda \lc \lc \ox'\shpr \ox''\rc\shpr \ox'''\rc
    +\lambda \lc \big(\ox'\shpr \ox''\big)\shpr \lc \ox'''\rc \rc
    + \lambda^2 \lc \big(\ox'\shpr \ox''\big) \shpr \ox'''\rc.
\end{eqnarray*}}
Similarly we obtain \allowdisplaybreaks{
\begin{eqnarray*} \frakx'
\shpr \big(\frakx''\shpr \frakx'''\big) &=& \lc\ox'\rc  \shpr
\Big(\lc\lc\ox''\rc \shpr \ox'''\rc
    + \lc\ox''\shpr \lc \ox'''\rc\rc +\lambda\lc \ox''\shpr\ox'''\rc \Big)\\
&=& \lc\lc\ox'\rc\shpr \big(\lc\ox''\rc\shpr\ox'''\big)\rc
    +\lc \ox'\shpr \lc \lc \ox''\rc \shpr \ox'''\rc\rc
    + \lambda \lc \ox'\shpr \big(\lc\ox''\rc\shpr \ox'''\big)\rc\\
&&  + \lc\lc \ox'\rc\shpr \big(\ox''\shpr \lc \ox'''\rc \big) \rc
    + \lc \ox' \shpr \lc \ox'' \shpr \lc \ox'''\rc\rc\rc
    + \lambda \lc \ox' \shpr \big(\ox''\shpr \lc \ox'''\rc \big)\rc\\
&&  + \lambda\lc\lc\ox'\rc\shpr \big( \ox''\shpr \ox'''\big) \rc
    + \lambda \lc \ox'\shpr \lc \ox''\shpr \ox'''\rc\rc
    + \lambda^2 \lc \ox' \shpr \big( \ox''\shpr\ox'''\big) \rc \Big)\\
&=& \lc\lc\lc\ox'\rc\shpr \ox''\rc\shpr \ox'''\rc
    + \lc \lc \ox'\shpr \lc\ox''\rc\rc \shpr \ox'''\rc
    + \lambda \lc\lc\ox'\shpr\ox''\rc\shpr\ox'''\rc\\
&&  +\lc \ox'\shpr \lc \lc \ox''\rc \shpr \ox'''\rc\rc
    + \lambda \lc \ox'\shpr \big(\lc\ox''\rc\shpr \ox'''\big)\rc\\
&&  + \lc\lc \ox'\rc\shpr \big(\ox''\shpr \lc \ox'''\rc \big) \rc
    + \lc \ox' \shpr \lc \ox'' \shpr \lc \ox'''\rc\rc\rc
    + \lambda \lc \ox' \shpr \big(\ox''\shpr \lc \ox'''\rc \big)\rc\\
&&  + \lambda\lc\lc\ox'\rc\shpr \big( \ox''\shpr \ox'''\big) \rc
    + \lambda \lc \ox'\shpr \lc \ox''\shpr \ox'''\rc\rc
    + \lambda^2 \lc \ox' \shpr \big( \ox''\shpr\ox'''\big) \rc.
\end{eqnarray*}}
Now by induction, the $i$-th term in the expansion of
$(\frakx'\shpr \frakx'')\shpr \frakx'''$ matches with the $\sigma(i)$-th
term  in the expansion of $\frakx'\shpr(\frakx'' \shpr \frakx''')$.
Here the permutation $\sigma\in \Sigma_{11}$ is
\begin{equation}
\left ( \begin{array}{c} i\\\sigma(i)\end{array}\right)
= \left ( \begin{array}{ccccccccccc} 1&2&3&4&5&6&7&8&9&10&11\\
    1&6&9&2&4&7&10&5&3&8&11\end{array} \right ).
\mlabel{eq:sigma}
\end{equation}
This completes the proof of the first part of Theorem~\mref{thm:freeao}.

\mref{it:RB}. The proof is immediate from the definition
$R_B(\frakx)=\lc \frakx\rc$ and Eq. (\mref{eq:shprod0}).

\mref{it:free}. Let $(A,R)$ be a unitary \rba of weight
$\lambda$. Let $f: B\to A$ be
a nonunitary $\bfk$-algebra morphism. We will construct a $\bfk$-linear map
$\free{f}:\ncsha(B)\to A$ by defining $\free{f}(\frakx)$ for
$\frakx\in \frakX_\infty$.
We achieve this by defining $\free{f}(\frakx)$ for $\frakx\in \frakX_n,\ n\geq 0$,
using induction on $n$.
For $\frakx\in \frakX_0:=X$, define
$\free{f}(\frakx)=f(\frakx).$ Suppose $\free{f}(\frakx)$ has been
defined for $\frakx\in \frakX_n$ and consider $\frakx$ in $\frakX_{n+1}$ which
is, by definition and Eq.~(\mref{eq:words2}), \allowdisplaybreaks{
\begin{eqnarray*} \altx(X,\frakX_{n})& =&
    \Big( \dbigcup_{r\geq 1} (X\lc \frakX_{n}\rc)^r \Big) \dbigcup
    \Big(\dbigcup_{r\geq 0} (X\lc \frakX_{n}\rc)^r  X\Big) \\
 &&    \dbigcup \Big( \dbigcup_{r\geq 0} \lc \frakX_{n}\rc (X\lc \frakX_{n}\rc)^r \Big)
   \dbigcup \Big( \dbigcup_{r\geq 0} \lc \frakX_{n}\rc (X\lc \frakX_{n}\rc)^r X\Big).
%\mlabel{eq:words2}
\end{eqnarray*}}
Let $\frakx$ be in the first union component
$\dbigcup_{r\geq 1} (X\lc \frakX_{n}\rc)^r$ above.
Then
$$\frakx = \prod_{i=1}^r(\frakx_{2i-1} \lc \frakx_{2i} \rc)$$
for
$\frakx_{2i-1}\in X$ and $\frakx_{2i}\in \frakX_n$, $1\leq i\leq r$.
By the construction of the multiplication $\shpr$ and the Rota-Baxter operator
$R_B$, we have
$$\frakx= \shpr_{i=1}^r(\frakx_{2i-1} \shpr \lc \frakx_{2i}\rc)
    = \shpr_{i=1}^r(\frakx_{2i-1} \shpr R_B(\frakx_{2i})).$$
Define
\begin{equation}
\free{f}(\frakx) = \ast_{i=1}^r \big(\free{f}(\frakx_{2i-1})
    \ast R\big (\free{f}(\frakx_{2i})) \big).
\mlabel{eq:hom}
\end{equation}
where the right hand side is well-defined by the induction hypothesis.
Similarly define $\free{f}(\frakx)$ if $\frakx$ is in the other union
components.
For any $\frakx\in \frakX_\infty$, we have
$R_B(\frakx)=\lc \frakx\rc\in \frakX_\infty$, and
by definition (Eq. (\mref{eq:hom})) of $\free{f}$, we have
\begin{equation}
\free{f}(\lc \frakx \rc)=R(\free{f}(\frakx)).
\mlabel{eq:hom1-2}
\end{equation}
So $\free{f}$ commutes with
the Rota-Baxter operators.
Combining this equation with Eq.~(\mref{eq:hom}) we see that if
$\frakx=\frakx_1\cdots \frakx_b$ is the standard decomposition of $\frakx$,
then
\begin{equation}
 \free{f}(\frakx)=\free{f}(\frakx_1)*\cdots * \free{f}(\frakx_b).
\mlabel{eq:staohom}
\end{equation}

Note that this is the only possible way to define $\free{f}(\frakx)$ in order
for $\free{f}$ to be a Rota-Baxter algebra homomorphism extending $f$.

We remain to prove that the map $\free{f}$ defined in Eq.~(\mref{eq:hom}) is
indeed an algebra homomorphism.
For this we only need to check the multiplicity
\begin{equation}
\free{f} (\frakx \shpr \frakx')=\free{f}(\frakx) \ast \free{f}(\frakx')
\mlabel{eq:hom2}
\end{equation}
for all $\frakx,\frakx'\in \frakX_\infty$.
For this we use induction on the sum of depths
$n:=d(\frakx)+d(\frakx')$. Then $n\geq 0$.
When $n=0$, we have $\frakx,\frakx'\in X$. Then Eq.~(\mref{eq:hom2}) follows from
the multiplicity of $f$. Assume the multiplicity holds for $\frakx,\frakx'
\in \frakX_\infty$ with $n\geq k$ and take $\frakx,\frakx'\in \frakX_\infty$ with
$n=k+1$.
Let $\frakx=\frakx_1\cdots \frakx_b$ and $\frakx'=\frakx'_1\cdots\frakx'_{b'}$
be the standard decompositions.
By Eq.~(\mref{eq:shprod0}),
\begin{align*}
\free{f}(\frakx_b\shpr \frakx'_1)&=
\left \{\begin{array}{ll}
\free{f}(\frakx_b \spr \frakx'_1), & {\rm if\ } \frakx_b,\frakx'_1\in X,\\
\free{f}(\frakx_b \frakx'_1), & {\rm if\ } \frakx_b\in X, \frakx'_1\in \lc \frakX_\infty\rc,\\
\free{f}(\frakx_b \frakx'_1), & {\rm if\ } \frakx_b\in \lc \frakX_\infty\rc,
    \frakx'_1\in X,\\
\free{f}\big( \lc \lc \ox_b\rc \shpr \ox'_1\rc +\lc \ox_b \shpr \lc \ox'_1\rc \rc
+\lambda \lc \ox_b \shpr \ox'_1 \rc\big), & {\rm if\ } \frakx_b=\lc \ox_b\rc,
\frakx'_1=\lc \ox'_1\rc \in \lc \frakX_\infty \rc.
\end{array} \right .
\end{align*}
In the first three cases, the right hand side is
$\free{f}(\frakx_b)*\free{f}(\frakx'_1)$ by the definition of $\free{f}$.
In the fourth case, we have, by Eq.~(\mref{eq:hom1-2}), the induction
hypothesis and the Rota-Baxter relation of $R$,
\begin{align*}
&\free{f}\big( \lc \lc \ox_b\rc \shpr \ox'_1\rc
    +\lc \ox_b \shpr \lc \ox'_1\rc \rc
+\lambda \lc \ox_b \shpr \ox'_1 \rc\big)\\
=&\free{f}(\lc \lc \ox_b\rc \shpr \ox'_1\rc)
+ \free{f}(\lc \ox_b \shpr \lc \ox'_1\rc \rc)
+\free{f}(\lambda \lc \ox_b \shpr \ox'_1 \rc)\\
=&R(\free{f}(\lc \ox_b\rc \shpr \ox'_1))
+ R(\free{f}(\ox_b \shpr \lc \ox'_1\rc ))
+ \lambda R(\free{f}(\ox_b \shpr \ox'_1 ))\\
=&R(\free{f}(\lc \ox_b\rc)*\free{f}(\ox'_1))
+ R(\free{f}(\ox_b) *\free{f}( \lc \ox'_1\rc ))
+ \lambda R(\free{f}(\ox_b) * \free{f}(\ox'_1) )\\
=&R(R(\free{f}(\ox_b))*\free{f}(\ox'_1))
+ R(\free{f}(\ox_b) *R(\free{f}(\ox'_1)))
+ \lambda R(\free{f}(\ox_b) * \free{f}(\ox'_1) )\\
=& R(\free{f}(\ox_b))*R(\free{f}(\ox'_1))\\
=& \free{f}(\lc \ox_b\rc) * \free{f}(\lc\ox'_1\rc)\\
=& \free{f} (\frakx_b) *\free{f}(\frakx'_1).
\end{align*}
Therefore $\free{f}(\frakx_b\shpr \frakx'_1)=\free{f}(\frakx_b)*\free{f}(\frakx'_1)$.
Then
\begin{align*}
\free{f}(\frakx\shpr \frakx')&=
\free{f}\big(\frakx_1\cdots\frakx_{b-1}(\frakx_b\shpr \frakx'_1)\frakx'_2\cdots
    \frakx'_{b'}\big) \\
&= \free{f}(\frakx_1)*\cdots *\free{f}(\frakx_{b-1})*
\free{f}(\frakx_b\shpr \frakx'_1)*\free{f}(\frakx'_2)\cdots
    \free{f}(\frakx'_{b'})\\
&= \free{f}(\frakx_1)*\cdots *\free{f}(\frakx_{b-1})*
\free{f}(\frakx_b)* \free{f} (\frakx'_1)*\free{f}(\frakx'_2)\cdots
    \free{f}(\frakx'_{b'})\\
&= \free{f}(\frakx)*\free{f}(\frakx').
\end{align*}
This is what we need.
\end{proof}

%%%%%%%%%%%%%%%%%%%%%%%%%%%%%%%%%%%%%%%%%%%%%%%%%%%%%%%%%%%%%%%%%%%%
%%%%%%%%%%%%%%%%%%%%%%%%%%%%%%%%%%%%%%%%%%%%%%%%%%%%%%%%%%%%%%%%%%%%%

\section{Universal enveloping algebras of dendriform trialgebras}
\mlabel{sec:adj}

\subsection{Dendriform dialgebras and trialgebras}
\mlabel{sec:dend}

We recall the following definitions.
A dendriform dialgebra~\mcite{Lo1} is a module $D$ with two binary operations
$\prec$ and $\succ$ such that
\begin{eqnarray}
&& (x \prec y) \prec z= x \prec (y\prec z +y \succ z),
(x \succ y ) \prec z= x \succ (y\prec z), \notag \\
&& (x \prec y +x\succ y)\succ z = x \succ (y\succ z)
\mlabel{eq:dia}
\end{eqnarray}
for $x,y,z\in D$.

A dendriform trialgebra~\mcite{L-R2} is
a module $T$ equipped with binary operations $\prec,\succ$ and
$\spr$ that satisfy the relations \allowdisplaybreaks{
\begin{eqnarray}
&&(x\prec y)\prec z=x\prec (y\star z),
(x\succ y)\prec z=x\succ (y\prec z),\notag \\
&&(x\star y)\succ z=x\succ (y\succ z),
(x\succ y)\spr z=x\succ (y\spr z),
\mlabel{eq:tri}\\
&&(x\prec y)\spr z=x\spr (y\succ z),
(x\spr y)\prec z=x\spr (y\prec z),
(x\spr y)\spr z=x\spr (y\spr z). \notag
\end{eqnarray}
Here
$\star=\prec+\succ+\spr.$
The category of dendriform trialgebras $(D,\prec,\succ,\spr)$ is denoted by $\DT$.
Recall that $\spr$, as well as $\star$, is an associative product.
The category $\Dend$ of dendriform dialgebras can be identified with the subcategory
of $\DT$ of objects with $\spr=0$.

These algebras are related to Rota-Baxter algebras by the following theorem.
\begin{theorem} {\bf (Aguiar~\mcite{Ag2}, Ebrahimi-Fard~\mcite{EF1})}
\begin{enumerate}
\item  A Rota-Baxter algebra $(A,R)$ of weight zero defines
a dendriform dialgebra $(A,\prec_R,\succ_R)$, where
\begin{equation}
 x\prec_R y=xR(y),\ x\succ_R y=R(x)y.
\mlabel{it:ags}
\end{equation}
\item A Rota-Baxter algebra $(A,R)$ of weight $\lambda$ defines a
dendriform trialgebra $(A,\prec_R,\succ_R,\spr_R)$, where
\begin{equation} x\prec_R y=xR(y),\ x\succ_R y=R(x)y, x\spr_R y=\lambda xy.
\mlabel{it:ef1s}
\end{equation}
\item  A Rota-Baxter
algebra $(A,R)$ of weight $\lambda$ defines a dendriform dialgebra
$(A,\prec'_R,\succ'_R)$, where
\begin{equation} x\prec'_R y=xR(y) + \lambda xy,\ x\succ'_R y=R(x)y.
\mlabel{it:ef2s}
\end{equation}
\end{enumerate}
\mlabel{thm:EFs}
\end{theorem}
We note that (\mref{it:ef2s}) specializes to (\mref{it:ags}) when $\lambda=0$.
The same can be said of (\mref{it:ef1s}) since when $\lambda=0$, the product
$\spr_R$ is zero and the relations of the trialgebra reduces to the relations
of a dialgebra.

It is easy to see that the maps between objects in the categories $\RBo_\lambda$,
$\Dend$ and $\DT$ in Theorem~\mref{thm:EFs} are
compatible with the morphisms.  Thus we obtain functors
$$\cale: \RBo_\lambda \to \DT,\ \calf: \RBo_\lambda \to \Dend.$$
We will study their adjoint functors. The two functors $\cale$ and
$\calf$ are related by the following simple observation.
\begin{prop}
\begin{enumerate}
\item
Let $(D,\prec,\succ,\spr)$ be in $\DT$. Then $(D,\prec',\succ')$ is in $\Dend$.
Here $\prec'=\prec+\spr$ and $\succ'=\succ$.
\mlabel{it:dte}
\item
Let $\calg:\DT \to \Dend$ be the functor obtained from \mref{it:dte}. Then
we have $\calf=\calg\circ \cale$.
\mlabel{it:dtf}
\item Fix a $\lambda\in \bfk$.
If the adjoint functors $\cale': \DT\to \RBo_\lambda$ and
$\calg':\Dend \to \DT$ exist, then the adjoint functor
$\calf':\Dend\to \RBo_\lambda$ exists and
$\calf'=\cale' \circ \calg'$.
\mlabel{it:dtc}
\end{enumerate}
\mlabel{pp:DtT}
\end{prop}
\begin{proof}
\mref{it:dte}
Let $\star' = \prec'+\succ$. Then we have $\star'=\star$. We have
\allowdisplaybreaks{
\begin{eqnarray*}
(a\prec' b)\prec' c &=& (a\spr b+a \prec b)\prec' c \\
&=& (a\spr b+a \prec b)\spr c + (a\spr b+a \prec b)\prec c \\
&=& (a\spr b) \spr c +(a \prec b)\spr c + (a\spr b)\prec c + (a \prec b)\prec c \\
&=& a\spr (b \spr c) +a \spr (b\succ c) + a \spr (b\prec c) + a \prec (b\star c)
\ \ {\rm (by\ Eq.~(\mref{eq:tri}))}\\
&=& a\prec' (b\star' c).
\end{eqnarray*}}
This verifies the first relation for the dendriform dialgebra.
The other two relations are also easy to verify:
$$ (a\succ' b)\succ' c = (a\succ b) \succ c= a \succ (b\star c)=a \succ' (b\star' c).$$
$$ (a\succ' b)\prec' c = (a\succ b)\spr c+(a\succ b)\prec c
=a\succ (b\spr c)+a \succ (b\prec c)=a \succ'(b\prec' c).$$

\mref{it:dtf}
For $(A,R)\in \RBo_\lambda$, by Theorem~\mref{thm:EFs} and item \mref{it:dte},
we have
\begin{eqnarray*}
\calg(\cale((A,R)))&=& \calg((A,\prec_R,\succ_R,\cdot_R))\\
&=& (A,\prec_R+\cdot_R, \succ_R)\\
&=& \calf((A,R)).
\end{eqnarray*}
It is easy to check that the composition is also compatible with the morphisms.
So we get the equality of functors.

\mref{it:dtc}
is standard: for any $C\in \Dend$ and $A\in \RBo_\lambda$, we have
\begin{eqnarray*}
\Hom(C,\calg(\calf(A))) &\cong & \Hom(\calg'(C),\calf(A))\\
&\cong& \Hom(\calf'(\calg'(C)),A).
\end{eqnarray*}
So $\calf'(\calg'(C))=\cale'(C)$.
\end{proof}

\subsection{Universal enveloping Rota-Baxter algebras}
\mlabel{ss:envel}

Motivated by the enveloping algebra of a Lie algebra, we are naturally led
to the following definition.
\begin{defn}
Let $D\in \DT$ (resp. $\Dend$) and let $\lambda\in \bfk$.
A {\bf universal enveloping Rota-Baxter algebra} of weight $\lambda$ of $D$ is
a Rota-Baxter algebra $\rbadj(D):=\rbadj_\lambda(D)\in \RBo_\lambda$ with a morphism
$\rho: D\to \rbadj(D)$
in $\DT$ (resp. $\Dend$) such that for any $A\in \RBo_\lambda$ and morphism
$f:D\to  A$ in $\DT$ (resp. $\Dend$), there is a unique
$\den{f}: \rbadj(D)\to A$ in $\RBo_\lambda$ such that $\den{f} \circ \rho =f$.
\mlabel{de:env}
\end{defn}
By the universal property of $\rbadj(D)$, it is unique up to isomorphisms in $\RBo_\lambda$.

\subsection{The existence of enveloping algebras}
We will separately consider the enveloping algebras for dialgebras
and trialgebras.

\subsubsection{The trialgebra case}
Let $D=(D,\prec,\succ,\spr)\in  \DT$. Then $(D,\spr)$ is a
nonunitary $\bfk$-algebra. Let $\lambda\in \bfk$
be given. Let $\ncshao(D):=\ncshao_\lambda(D)$ be the free nonunitary
Rota-Baxter algebra over $D$ of weight $\lambda$ constructed in
\S\mref{ss:prodao}. Identify $D$ as a subalgebra of $\ncshao(D)$
by the natural injection $j_D$ in Eq.(\mref{eq:jo}). Let $I_R$ be
the Rota-Baxter ideal of $\ncshao(D)$ generated
by the set
\begin{equation}
\big \{ x\prec y - x\lc y\rc,\; x\succ y - \lc x\rc y\  \big|\
x,y\in D \big\}. \mlabel{eq:gen}
\end{equation}
Here a Rota-Baxter ideal of $\ncshao(D)$ is an ideal $I$ of $\ncshao(D)$ such that
$R_B(I)\subseteq I$, and the Rota-Baxter ideal of $\ncshao(D)$ generated by
a subset of $\ncshao(D)$ is the intersection of all Rota-Baxter ideals of $\ncshao(D)$ that
contain the subset.
Let $\pi: \ncshao(D)\to \ncshao(D)/I_R$ be the quotient map.

\begin{theorem}
The quotient Rota-Baxter algebra $\ncshao(D)/I_R$, together with
$\rho:=\pi \circ j_D$, is the universal enveloping Rota-Baxter algebra of $D$.
\mlabel{thm:env}
\end{theorem}
The theorem provides the adjoint functor $\cale':\DT \to \RBo$ of the
functor $\cale: \RBo\to \DT$.
\begin{proof}
Let $(A,R)\in \RBo_\lambda$. It gives an object in $\DT$ by Theorem~\mref{thm:EFs}
which we still denote by $A$. Let $f:D\to A$ be a morphism
in $\DT$. We will complete the following commutative diagram
\begin{equation}
\xymatrix{
D \ar[rr]^{j_D} \ar[d]_f  && \ncshao(D) \ar[d]^\pi \ar@{.>}[dll]_{\free{f}} \\
A && \ncshao(D)/I_R \ar@{.>}[ll]_{\den{f}}
}
\end{equation}
By the freeness of $\ncshao(D)$, there is a morphism
$\free{f}:\ncshao(D) \to A$ in $\RB^0$ such that the upper left
triangle commutes. So for any $x,y\in D$, by Eq. (\mref{eq:hom}), we have
\allowdisplaybreaks{
\begin{eqnarray*}
\free{f}(x\prec y - x\lc y \rc)
&=& \free{f}(x\prec y) - \free{f}(x\lc y\rc) \\
&=&\free{f}(x\prec y)-\free{f}(x)R(\free{f}(y))\\
 &=& f(x\prec y) -f(x)R(f(y))\\
&=& f(x\prec y)-f(x)\prec_R f(y)\\
&=& f(x\prec y)-f(x\prec y)=0.
\end{eqnarray*}}
Therefore, $x\prec y - x\lc y\rc$ is in $\ker(\free{f})$. Similarly,
$x\succ y -\lc x\rc y$ is in $\ker(\free{f})$. Thus $I_R$ is in $\ker(\free{f})$
and there is a morphism $\den{f}: \ncshao(D)/I_R\to A$ in $\RBo$ such that
$\free{f}=\den{f} \circ \pi$.
Then
$$ \den{f}\circ \rho = \den{f} \circ \pi \circ j_D=\free{f}\circ j_D=f.$$
This proves the existence of $\den{f}$.

Suppose $\den{f}':\ncshao(D)/I_R \to A$ is a morphism in $\RBo$ such that
$\den{f}'\circ \rho=f$. Then
$$ (\den{f}' \circ \pi)\circ j_D = f = (\den{f}\circ \pi)\circ j_D.$$
By the universal property of the free Rota-Baxter algebra $\ncshao(D)$ over $D$,
we have $\den{f}'\circ \pi = \den{f} \circ \pi$ in $\RBo$. Since $\pi$ is
surjective, we have $\den{f}'=\den{f}$. This proves the uniqueness of $\den{f}$.
\end{proof}

\subsubsection{The dialgebra case}
Now let $D=(D,\prec,\succ)\in \Dend$.
Let $T(D)=\bigoplus_{n\geq 1} D^{\ot n}$ be the tensor product algebra over $D$.
Then $T(D)$ is the free nonunitary algebra generated by the $\bfk$-module $D$~\cite[Prop. II.5.1]{Ka}.
By Corollary~\mref{co:vecfree},
$\ncshao(T(D))$, with the natural injection
$i_D: D\to T(D) \to \ncshao(T(D))$, is the free Rota-Baxter
algebra over the vector space $D$.

Let $J_R$ be the Rota-Baxter ideal of $\ncshao(T(D))$ generated
by the set
\begin{equation}
\big \{ x\prec y - x\lc y\rc-\lambda x\ot y,\;
    x\succ y - \lc x\rc y\  \big|\ x,y\in D \big\}
\mlabel{eq:gendend}
\end{equation}
Let $\pi: \ncshao(T(D))\to \ncshao(T(D))/J_R$ be the quotient map.

\begin{theorem}
The quotient Rota-Baxter algebra $\ncshao(T(D))/J_R$, together
with  $\rho:= \pi \circ i_D$, is the universal
enveloping Rota-Baxter algebra of $D$ of weight $\lambda$.
\mlabel{thm:envdend}
\end{theorem}

\begin{proof}
Let $(A,R)$ be a Rota-Baxter algebra of weight $\lambda$ and let $f:D\to A$ be
a morphism in $\Dend$. More precisely, we have $f:D\to \calg A$ where
$\calg A=(A,\prec_R',\succ_R')$ is the dendriform dialgebra in Theorem~\mref{thm:EFs}.
We will complete the following commutative diagram, using notations from
Corollary~\mref{co:vecfree}.
\begin{equation}
\xymatrix{ & T(D) \ar[rd]^{j_{T(D)}} \ar@{.>}[lddd]^{\freea{f}} & \\
D \ar[rr]^{i_D} \ar[dd]_f \ar[ru]^{k_D} && \ncshao(T(D)) \ar[dd]^\pi \ar@{.>}[ddll]_{\freev{f}} \\
&& \\
A && \ncshao(T(D))/J_R \ar@{.>}[ll]_{\den{f}}
}
\end{equation}

By the universal property of the free algebra $T(D)$ over $D$, there is a
unique morphism $\freea{f}:T(D)\to A$ in $\Algo$ such that
$\freea{f}\circ k_D =f$ and so
$\freea{f}(x_1\ot \cdots \ot x_n)=f(x_1) * \cdots * f(x_n)$.
Here $*$ is the product in $A$.
Then by the universal property of the free Rota-Baxter algebra
$\ncshao(T(D))$ over $T(D)$, there is a unique morphism
$\free{\freea{f}}:\ncshao(T(D)) \to A$ in $\RBo$
such that $\free{\freea{f}}\circ j_{T(D)} =\freea{f}$.
By Corollary~\mref{co:vecfree}, $\free{\freea{f}}=\freev{f}.$
Then
\begin{equation} \freev{f}\circ i_D =\freev{f} \circ j_{T(D)} \circ k_D
= \freea{f} \circ k_D = f.
\mlabel{eq:free2}
\end{equation}
So
for any $x,y\in D$, we have
\begin{eqnarray*}
\freev{f}(x\prec y - x\lc y \rc-\lambda x\ot y)
&=&\freev{f}(x\prec y)-\freev{f}(x)*R(\freev{f}(y))-\lambda \freev{f}(x\ot y)\\
&=&\freev{f}(x\prec y)-\freev{f}(x)*R(\freev{f}(y))-\lambda \freea{f}(x\ot y)\\
&=& f(x\prec y) -f(x)*R(f(y)) -\lambda f(x)* f(y)\\
&=& f(x\prec y)-f(x)\prec_R' f(y)\\
&=& f(x\prec y)-f(x\prec y)=0.
\end{eqnarray*}
Therefore, $x\prec y - x\lc y\rc-\lambda x\ot y$ is in $\ker(\freev{f})$. Similarly,
$x\succ y -\lc x\rc y$ is in $\ker(\freev{f})$. Thus $J_R$ is in $\ker(\freev{f})$
and there is a morphism $\den{f}: \ncshao(T(D))/J_R\to A$ in $\RBo$ such that
$\freev{f}=\den{f} \circ \pi$.
Then by the definition of $\rho=\pi \circ i_D$ in the theorem and
Eq. (\mref{eq:free2}), we have
$$ \den{f}\circ \rho = \den{f} \circ \pi \circ i_D=\freev{f}\circ i_D=f.$$
This proves the existence of $\den{f}$.

Suppose $\den{f}':\ncshao(T(D))/J_R \to A$ is also a morphism in $\RBo$
such that $\den{f}'\circ \rho=f$. Then
$$ (\den{f}' \circ \pi)\circ i_D = f = (\den{f}\circ \pi)\circ i_D.$$
By Corollary~\mref{co:vecfree}, the free Rota-Baxter algebra $\ncshao(T(D))$
over the algebra $T(D)$ is also the free Rota-Baxter algebra over the
vector space $D$ with respect the natural injection $i_D$.
So we have $\den{f}'\circ \pi = \den{f} \circ \pi$ in $\RBo$. Since $\pi$ is
surjective, we have $\den{f}'=\den{f}$. This proves the uniqueness of $\den{f}$.
\end{proof}

\section{Free dendriform di- and trialgebras and free Rota-Baxter algebras}
\mlabel{sec:dfree}
The results in this section can be regarded as more precise forms of results
in \S\mref{sec:adj} in special cases.
Our emphasis here is to interpret free dendriform dialgebras
and free dendriform trialgebras as natural subalgebras of free Rota-Baxter
algebras. This interpretation also suggests a planar tree structure on free Rota-Baxter algebras which will be made precise in~\mcite{free}.

\subsection{The dialgebra case}
\subsubsection{Free dendriform dialgebras}

Let $\bfk$ be a field.
We briefly recall the construction of free dendriform dialgebra $\Dend(V)$
over a $\bfk$-vector space $V$ as colored planar binary trees. For details, see
\mcite{Lo1,Ron}.

Let $X$ be a basis of $V$. For $n\geq 0$, let $Y_n$ be the set of
planar binary trees with $n+1$ leaves and one root such that the
valence of each internal vertex is exactly two. Let $Y_{n,X}$ be
the set of planar binary trees with $n+1$ leaves and with vertices
decorated by elements of $X$. The unique tree with one leave is
denoted by $|$. So we have $Y_0=Y_{0,X}=\{|\}$. Let
$\bfk[Y_{n,X}]$ be the $\bfk$-vector space generated by $Y_{n,X}$.
Here are the first few of them without decoration.
$$Y_0 = \{ \ \vert\ \} ,\qquad \ Y_1 = \Big\{\ \begin{array}{c} \\[-0.4cm] \shu1 \\ \end{array}\Big\} ,
\qquad
Y_2 = \Big\{\ \begin{array}{c} \\[-0.4cm] \shub2\ , \shuc3\ \\ \end{array} \Big\}
$$
\allowdisplaybreaks{\begin{eqnarray*}
Y_3 &=& \Big\{\ \begin{array}{c} \\[-0.4cm] \5tree,\
\6tree,\ \c1tree,\ \ldots
\\ \end{array} \Big\}.
\end{eqnarray*}}
%
%Here are the first few of them without decoration.
%$$ \Big\vert\  ,\qquad \  \arbreA, \qquad  \arbreB , \qquad \arbreC, \qquad
%\arbretrois $$
%
For $T\in Y_{m,X},U\in Y_{n,X}$ and $x\in X$, the grafting of $T$ and $U$ over
$x$ is $T\vee_x U\in Y_{m+n+1,X}$.
Let $\Dend(V)$ be the graded vector space $\bigoplus_{n\geq 1} \bfk[Y_{n,X}]$.
Define binary operations $\prec$ and $\succ$ on $\Dend(V)$ recursively by
\begin{enumerate}
\item
$|\succ T = T\prec |=T$ and $|\prec T = T\succ |=0$ for $T\in Y_{n,X}, n\geq 1$;
\item
For $T=T^\ell\vee_x T^r$ and $U=U^\ell\vee_y U^r$, define
$$ T\prec U= T^\ell \vee_x (T^r\prec U+T^r \succ U),\quad
  T\succ U = (T\prec U^\ell+T\succ U^\ell) \vee_y U^r.$$
\end{enumerate}
Since $|\prec |$ and $|\succ |$ is not defined, the binary operations $\prec$ and
$\succ$ are only defined on $\Dend(V)$ though the operation $\star:=\prec +\succ$
can be extended to $H_\LR:=\bfk[Y_0]\oplus \Dend(V)$ by defining $|\star T=T\star |=T.$
By~\mcite{Lo1}} $(\Dend(V),\prec,\succ)$ is the free
dendriform dialgebra over $V$.

\begin{theorem}
Let $V$ be a $\bfk$-vector space. The free dendriform dialgebra over
$V$ is a sub dendriform dialgebra of the free Rota-Baxter algebra
$\ncshao(V)$ of weight zero.
\mlabel{thm:dial}
\end{theorem}

The proof will be given in the next subsection.

\subsubsection{Proof of Theorem~\mref{thm:dial}}

For the given vector space $V$, make $V$ into a $\bfk$-algebra without identity
by given
$V$ the zero product. Let $\ncshao(V)$ be the free nonunitary Rota-Baxter algebra
of weight zero over
$V$ constructed in Theorem~\ref{thm:freeao}. Since $\ncshao(V)$
is a dendriform dialgebra, the natural map $j_V: V\to \ncshao(V)$ extends
uniquely to
a dendriform dialgebra morphism $D(j): \Dend(V)\to \ncshao(V)$.
We will prove that this map is injective and identifies $\Dend(V)$ as a subalgebra
of $\ncshao(V)$ in the category of dendriform dialgebras.
We first define a map
$$\phi: \Dend(V)\to \ncshao(V)$$
and then show in Theorem~\mref{thm:freedend} below that it agrees with
$D(j)$.
We construct $\phi$ by defining $\phi(T)$ for $T\in Y_{n,X}, n\geq 1,$ inductively
on $n$. Any $T\in Y_{n,X}, n\geq 1$ can be uniquely written as
$T=T^\ell \vee_x T^r$ with $x\in X$ and $T^\ell,T^r\in \cup_{0\leq i<n} Y_{i,X}$.
We then define
\begin{equation}
\phi(T)=\left \{\begin{array}{ll}
\lc \phi(T^\ell)\rc x \lc \phi(T^r) \rc, & T^\ell\neq 1, T^r \neq 1,\\
x \lc \phi(T^r) \rc, &  T^\ell =1, T^r\neq 1,\\
\lc \phi(T^\ell) \rc x, & T^\ell\neq 1, T^r =1,\\
x, & T^\ell=1,T^r=1.
\end{array} \right .
\end{equation}
For example,
$$
%\phi\Big ( \arbreAd \Big)
\phi \Big(\ \begin{array}{c} \\[-0.4cm] \xtree \\ \end{array} \Big)= x,
\qquad \phi \Bigg(\ \begin{array}{c} \\[-0.4cm] \xyztree \\ \end{array}\Bigg)
%\arbretroisd
= \lc x \rc z \lc y \rc.
$$ We recall~\mcite{Lo1} that $\Dend(V)$ with the operation
$\star:=\,\prec+\succ$ is an associative algebra.
%We also note that, for any Rota-Baxter algebra $(A,R)$
%of weight 0, the new product $a*_R b:=aR(b)+R(a)b=a\prec b+a\succ b$ makes $A$
%into an associative algebra (without a unit). Further we have
%$R(a)R(b)=R(a*_Rb)$.

We now describe a submodule of $\ncshao(V)$ to be identified with the
image of $\phi$ in Theorem~\mref{thm:freedend}.

\begin{defn} {\rm
A $y\in \frakX_\infty$ is called a
{\bf dendriform diword (DW)} if it satisfies the following {\em
additional} properties.
\begin{enumerate}
%\item There is no subword $\lc \bfone \rc$ in the word;
\item $y$ is not in $\lc \frakX_\infty \rc$;
\item There is no subword $\lc \lc \frakx \rc \rc$ with $\frakx\in \frakX_\infty$
in the word;
\item There is no subword of the form
$\frakx_1\lc \frakx_2 \rc \frakx_3$ with $\frakx_1,\frakx_3\in X$ and
$\frakx_2\in \frakX_\infty$.
\end{enumerate}
We let $DW(V)$ be the subspace of $\ncshao(V)$ generated by the
dendriform diwords. }
%The set of dendriform diwords is denoted by $DW(V)$.
\end{defn}

For example
$$ x_0\lc x_1\lc x_2 \rc \rc,
    \lc x_0\rc x_1\lc x_2 \rc$$
are dendriform diwords while
$$
%\lc \bfone \rc,
\lc \lc x_1\rc \rc, \lc \lc x_1\rc x_2 \lc x_3\rc\rc,
     x_1 \lc x_2 \rc x_3  $$
are in $\frakX_\infty$ but not dendriform diwords.

Equivalently, $DW(V)$ can be characterized in terms of the
decomposition (\mref{eq:words3}). For subsets $Y,Z$ of $\frakX_\infty$, define
$$ D(Y,Z)=(Y\lc Z \rc ) \bigcup (\lc Z\rc Y) \bigcup \lc Z\rc Y \lc Z\rc.$$
Then define $D_0(V)=X$ and, for $n\geq 0$, inductively define
\begin{equation}
 D_{n+1}(V)=D(X,D_n(V))
= (X\lc D_n(V) \rc ) \bigcup (\lc D_n(V) \rc X) \bigcup \lc D_n(V)\rc X \lc D_n(V)\rc.
\mlabel{eq:dend}
\end{equation}
Then  $D_\infty: = \cup_{n\geq 0} D_n(V)$ is the set of dendriform
diwords and $DW(V)=\oplus_{\frakx\in D_\infty} \bfk \frakx.$

Theorem~\mref{thm:dial} follows from the following theorem.
\begin{theorem}
\begin{enumerate}
\item $\phi: \Dend(V)\to \ncshao(V)$ is a homomorphism of
dendriform dialgebras. \mlabel{it:hom}
\item $\phi=D(j)$, the morphism of dendriform dialgebras induced
by $j:V\to \ncshao(V)$. \mlabel{it:agree}
\item $\phi(\Dend(V))=DW(V)$.
\mlabel{it:image}
\item $\phi$ is injective. \mlabel{it:injd}
\end{enumerate}
\mlabel{thm:freedend}
\end{theorem}

\begin{proof}
\mref{it:hom}
we first note that the operations $\prec$ and $\succ$ can be equivalently defined
as follows.
Let $T\in Y_{m,X}, U\in Y_{n,X}$ with $m\geq 1, n\geq 1$. Then
$T=T^\ell \vee_x T^r, U=U^\ell \vee_y U^r$ with $x,y\in X$ and
$T^\ell,T^r,U^\ell,U^r\in \cup_{i\geq 0} Y_{i,X}.$ Define
\begin{eqnarray}
T\prec U: &=& \left \{ \begin{array}{ll}
    T^\ell \vee_x (T^r \prec U + T^r \succ U), &{\rm if\ } T^r\neq |,\\
    T^\ell \vee_x U,& {\rm if\ } T^r=|.
    \end{array} \right .
\\
T\succ U: &=&\left \{ \begin{array}{ll}
    (T\prec U^\ell+T\succ U^\ell) \vee_y U^r, & {\rm if\ } U^\ell \neq |,\\
    T \vee_y U^r, & {\rm if\ } U^\ell = |.
    \end{array} \right .
\end{eqnarray}
Thus we have
\begin{eqnarray*}
\phi(T\prec U) &=& \left \{ \begin{array}{ll}
    \phi( T^\ell \vee_x (T^r \prec U + T^r \succ U)), &{\rm if\ } T^r\neq |,\\
    \phi(T^\ell \vee_x U), & {\rm if\ } T^r=|.
    \end{array} \right . \\
&=& \left \{\begin{array}{ll}
 \lc \phi(T^\ell)\rc x \lc \phi(T^r \prec U + T^r \succ U)\rc,
                &{\rm if\ } T^r\neq |, T^\ell\neq |,\\
  x \lc \phi(T^r \prec U + T^r \succ U)\rc,
                &{\rm if\ } T^r\neq |, T^\ell = |,\\
\lc \phi(T^\ell)\rc x \lc \phi(U) \rc, &{\rm if\ } T^r = |, T^\ell\neq |,\\
 x \lc \phi(U)\rc, &{\rm if\ } T^r = |, T^\ell = |.
\end{array} \right . \\
&&{\rm(by\ definition\ of\ }\phi{)} \\
&=& \left \{\begin{array}{ll}
 \lc \phi(T^\ell)\rc x \lc \phi(T^r) \prec_R \phi(U) + \phi(T^r) \succ_R \phi(U)\rc,
                &{\rm if\ } T^r\neq |, T^\ell\neq |,\\
  x \lc (\phi(T^r) \prec_R \phi(U) + \phi(T^r) \succ_R \phi(U))\rc,
                &{\rm if\ } T^r\neq |, T^\ell = |,\\
\lc \phi(T^\ell)\rc x \lc \phi(U)\rc , &{\rm if\ } T^r = |, T^\ell\neq |,\\
 x \lc \phi(U)\rc, &{\rm if\ } T^r = |, T^\ell = |.
\end{array} \right . \\
&&{\rm(by\ induction\ hypothesis)}
\end{eqnarray*}

On the other hand, we have
\begin{eqnarray*}
\phi(T) \prec_R \phi(U)&=& \phi(T^\ell\vee_x T^r) \lc \phi(U) \rc  \\
&=& \left \{\begin{array}{ll}
    \lc \phi(T^\ell) \rc x \lc \phi(T^r)\rc\lc \phi(U)\rc,
            &{\rm\ if\ }T^r\neq |, T^\ell \neq |, \\
    x \lc \phi(T^r)\rc \lc \phi(U)\rc, &{\rm if\ } T^r\neq |, T^\ell=|, \\
    \lc \phi(T^\ell) \rc x \lc \phi(U) \rc ,&{\rm if\ } T^r=|, T^\ell\neq |,\\
    x\lc \phi(U) \rc, & {\rm if\ } T^r=|,T^\ell=|.
\end{array} \right . \\
&& {\rm (by\ definition\ of\ }\phi{)} \\
&=& \left \{\begin{array}{ll}
    \lc \phi(T^\ell) \rc x
    \big \lc \phi(T^r)\lc \phi(U)\rc +\lc \phi(T^r) \rc \phi(U) \big\rc,
            &{\rm\ if\ }T^r\neq |, T^\ell \neq |, \\
    x \lc \phi(T^r) \big \lc \phi(U)\rc +\lc \phi(T^r) \rc \phi(U) \big \rc,
        &{\rm if\ } T^r\neq |, T^\ell=|, \\
    \lc \phi(T^\ell) \rc x \lc \phi(U) \rc ,&{\rm if\ } T^r=|, T^\ell\neq |,\\
    x\lc \phi(U) \rc, & {\rm if\ } T^r=|,T^\ell=|.
\end{array} \right . \\
&& {\rm (by\ Rota-Baxter\ relation\ of\ }R(T)=\lc T\rc{)}. \\
\end{eqnarray*}
This proves $\phi(T\prec U)=\phi(T)\prec_R \phi(U)$. We similarly prove
$\phi(T\succ U) =\phi(T)\succ_R \phi(U).$
Thus $\phi$ is a homomorphism in $\Dend$.

\mref{it:agree} follows from the uniqueness of the dendriform dialgebra
morphism $\Dend(V)\to \ncshao(V)$ extending the map
$j_V:V\to \ncshao(V)$.

\mref{it:image} We only need to prove $DW(V)\subseteq  \phi(\Dend(V))$
and $\phi(\Dend(V)) \subseteq DW(V)$. To prove the former, we prove
$D_n\subseteq \phi(\Dend(V))$ by induction on $n$.

When $n=0$, $D_n=X$ so the inclusion is clear.
Suppose the inclusion holds for $n$. Then by the definition of $D_{n+1}(V)$
in Eq.~(\mref{eq:dend}), an element of $D_{n+1}(V)$ is of the following three
forms:

i) It is $\frakx \lc \frakx'\rc$ with $\frakx\in X$, $\frakx'\in D_n(V)$.
Then  it is $\frakx \prec_R \frakx'$ which is in $\phi(\Dend(V))$ by the induction
hypothesis and the fact that $\phi(\Dend(V))$ is a sub dendriform algebra.

ii) It is $\lc \frakx \rc \frakx'$ with $\frakx\in D_n(V)$ and $\frakx'\in X$.
Then the same proof works.

iii) It is $\lc \frakx\rc \frakx' \lc \frakx''\rc$ with $\frakx,\frakx''\in D_n(V)$
and $\frakx'\in X$. Then it is
$$(\frakx \succ_R \frakx') \prec_R \frakx''=\frakx' \succ_R (\frakx' \prec_R \frakx'').$$
By induction, $\frakx$ and $\frakx''$ are in the sub dendriform dialgebra
$\phi(\Dend(V))$. So the element itself is in $\phi(\Dend(V))$.

The second inclusion follows easily by induction on degrees of trees in
$\Dend(V)$.

\mref{it:injd} By the definition of $\phi$ and part \mref{it:image}, $\phi$ gives
a one-one correspondence between
$\cup_{n\geq 0} Y_{n,X}$ as a basis of $\Dend(V)$ and $DW(V)$ as a basis of
$\phi(\Dend(V))$. Therefore $\phi$ is injective.
\end{proof}

%%%%%%%%%%%%%%%%%%%%%%%%%%%%%%%%%%%%%%%%%%%%
%%%%%%%%% Trialgebras
%%%%%%%%%%%%%%%%%%%%%%%%%%%%%%%%%%%%%%%%%%%%%%%%%%%%%%%%%%%%%%

\subsection{The trialgebra case}
\subsubsection{Free dendriform trialgebras}

We describe the construction of free dendriform trialgebra $\DT(V)$
over a vector space $V$ as colored planar trees. For details when $V$ is
of rank one over $\bfk$,  see~\mcite{L-R1}.

Let $\Omega$ be a basis of $V$. For $n\geq 0$, let $T_n$ be the set of
planar trees with $n+1$ leaves and one root such that the valence
of each internal vertex is at least two. Let $T_{n,\Omega}$ be the set
of planar  trees with $n+1$ leaves and
with vertices {\bf valently decorated} by elements of $\Omega$, in the
sense that if a vertex has valence $k$, then the vertex is
decorated by a vector in $\Omega^{k-1}$. For example the vertex of
%$\arbreA$
\sshu1 is  decorated by $x\in \Omega$ while the vertex of
%$\arbreBC$
\sc4tree is decorated by $(x,y)\in \Omega^2.$ The unique tree with one
leaf is denoted by $|$. So we have $T_0=T_{0,\Omega}=\{|\}$. Let
$\bfk[T_{n,\Omega}]$ be the $\bfk$-vector space generated by $T_{n,\Omega}$.

Here are the first few of them without decoration.
$$T_0 = \{ \ \vert\ \} ,\qquad \ T_1 = \Big\{\ \begin{array}{c} \\[-0.4cm] \shu1 \\ \end{array}\Big\} ,
\qquad
T_2 = \Big\{\ \begin{array}{c} \\[-0.4cm] \shub2\ , \shuc3\ , \shud4\ \\ \end{array} \Big\}
$$
\allowdisplaybreaks{\begin{eqnarray*}
T_3 &=& \Big\{\ \begin{array}{c} \\[-0.4cm] \5tree,\
\6tree,\ \7tree,\ \8tree,\ \9tree,\ \cdots,\ \b1tree,\ \c1tree,\
\d1tree,\ \e1tree,\ \ldots \\ \end{array} \Big\}.
\end{eqnarray*}}

%
%Here are the first few of them without decoration.
%$$T_0 = \{\ \vert\ \} ,\qquad \ T_1 = \Big\{\ \arbreA\ \Big\} ,\qquad
%T_2 = \Big\{\ \arbreB , \arbreC ; \arbreBC\ \Big\}$$
%%\begin{eqnarray*}
%%T_3 &=& \Big\{\ \arbreun, \arbredeux, \arbretrois, \arbrequatre, \arbrecinq; \\
%%&&\arbreuut, \arbretut, \arbretuu, \arbreutt, \arbrettu; \arbrettt \Big\}.
%%\end{eqnarray*}
%

For $T^{(i)}\in T_{n_i,\Omega},\, 0\leq i\leq k,$ and $x_i\in \Omega,\, 1\leq i\leq k$,
the grafting of $T^{(i)}$ over $(x_1,\cdots,x_k)$ is
$$T^{(0)}\vee_{x_1} T^{(1)}\vee_{x_2} \cdots \vee_{x_k} T^{(k)}.$$
Any tree can be uniquely expressed as such a grafting of lower degree trees.
For example
$$ \dec1tree \ \ \ \ = | \vee_x | \vee_y |.$$
%$$ \arbreBCd \ \ \ \ = | \vee_x | \vee_y |.$$
Let $\DT(V)$ be the graded vector space $\bigoplus_{n\geq 1} \bfk[T_{n,\Omega}]$.
Define binary operations $\prec$, $\succ$ and $\spr$ on $\DT(V)$ recursively by
\begin{enumerate}
\item
$|\succ T = T\prec |=T$, $|\prec T = T\succ |=0$ and
$| \spr T = T\spr |=0$ for $T\in T_{n,\Omega}, n\geq 1$;
\item
For $T=T^{(0)}\vee_{x_1}\cdots \vee_{x_m} T^{(m)}$ and
$U=U^{(0)}\vee_{y_1}\cdots \vee_{y_n} U^{(n)}$, define
\begin{eqnarray*}
T\prec U &=& T^{(0)}\vee_{x_1} \cdots \vee_{x_m} (T^{(m)} \star U),\\
T\succ U &=& (T \star U^{(0)}) \vee_{y_1}\cdots \vee_{y_n} U^{(n)},\\
T\spr U &=&  T^{(0)}\vee_{x_1} \cdots \vee_{x_m} (T^{(m)} \star U^{(0)})
    \vee_{y_1}\cdots \vee_{y_n} U^{(n)}.
\end{eqnarray*}
\end{enumerate}
Here $\star:=\prec +\succ+\, \spr$
Since $|\prec |$, $|\succ |$ and $|\spr |$ are not defined, the binary operations
$\prec$, $\succ$ and $\spr$ are only defined on $\DT(V)$ though the operation $\star$
can be extended to $H_\DT:=\bfk[T_0]\oplus \DT(V)$ by defining $|\star T=T\star |=T.$

\begin{theorem} $(\DT(V),\prec,\succ,\spr)$ is the free dendriform trialgebra
over $V$.
\mlabel{thm:LR}
\end{theorem}
\begin{proof}
The proof is given by Loday and Ronco in~\mcite{L-R1} when $V$ is of dimension one.
The proof for the general case is the same.
\end{proof}

Our goal is to prove
\begin{theorem}
Let $V$ be a $\bfk$-vector space. The free dendriform trialgebra over
$V$ is a canonical sub-dendriform trialgebra of the free Rota-Baxter algebra
$\ncshao(T(V))$ of weight one.
\mlabel{thm:tri}
\end{theorem}
We restrict the weight of the Rota-Baxter algebra to one to ease the notations.
The proof will be given in the next subsection.

\subsubsection{Proof of Theorem~\mref{thm:tri}}
Let $V$ be the given $\bfk$-vector space with basis $\Omega$.
Let $T(V)=\bigoplus_{n\geq 1} V^{\ot n}$ be the tensor product algebra over $V$.
Then $T(V)$ is the free nonunitary algebra generated by the $\bfk$-space $V$.
A basis of $T(V)$ is $X:=M(\Omega)$, the free semigroup generated
by $\Omega$.
By Theorem~\mref{thm:freeao}, $\ncshao(T(V)):=\ncshao_\bfone
(T(V))$ is the free nonunitary Rota-Baxter algebra over $T(V)$ of
weight $\bfone$ constructed in \S\mref{ss:prodao}.

Since $\ncshao(T(V))$
is a dendriform trialgebra, the natural map $j_V: V\to \ncshao(T(V))$ extends
uniquely to
a dendriform trialgebra morphism $T(j): \DT(V)\to \ncshao(T(V))$.
We will prove that this map is injective and identifies $\DT(V)$ as a subalgebra
of $\ncshao(T(V))$ in the category of dendriform trialgebras.
We first define a map
$$\psi: \DT(V)\to \ncshao(T(V))$$
and then show in Theorem~\mref{thm:freetri} below that it agrees with
$T(j)$.
We construct $\psi$ by defining $\psi(T)$ for $T\in T_{n,\Omega}, n\geq 1,$ inductively
on $n$.
Any $T\in T_{n,\Omega}, n\geq 1$, can be uniquely written as
$T=T^{(0)} \vee_{x_1} \cdots \vee_{x_k} T^{(k)}$ with $x_i\in \Omega$ and
$T^{(i)}\in \cup_{0\leq i<n} T_{i,\Omega}$.
We then define
\begin{equation}
\psi(T)=
\overline{\lc\psi(T^{(0)})\rc} x_1 \overline{\lc\psi(T^{(1)})\rc}
    \cdots \overline{\lc\psi(T^{(k-1)})\rc} x_k \overline{\lc \psi(T^{(k)})\rc},
\mlabel{eq:psi}
\end{equation}
where $\overline{\lc\psi(T^{(i)})\rc}=\lc\psi(T^{(i)})\rc$ if
$\psi(T^{(i)}) \neq |$. If $\psi(T^{(i)}) = |$, then the factor
$\lc\psi(T^{(i)})\rc$ is dropped when $i=0$ or $k$, and is replaced by
$\otimes$ when $0<i<k$. For example,
$$
\overline{\lc\psi(|)\rc} x_1 \overline{\lc\psi(T^{(1)})\rc} x_2
    \cdots  x_k \overline{\lc \psi(T^{(k)})\rc}
= x_1 \overline{\lc\psi(T^{(1)})\rc} x_2
    \cdots  x_k \overline{\lc \psi(T^{(k)})\rc}$$
and
$$
\overline{\lc\psi(T^{(0)})\rc} x_1 \overline{\lc\psi(|)\rc} x_2
    \overline{\lc\psi(T^{(2)})\rc} \cdots  x_k \overline{\lc \psi(T^{(k)})\rc}
= \overline{\lc\psi(T^{(0)})\rc} (x_1 \otimes x_2)
    \overline{\lc\psi(T^{(2)})\rc}\cdots  x_k \overline{\lc \psi(T^{(k)})\rc}.$$
In particular,
$$ \psi \bigg(\ \begin{array}{c} \\[-0.4cm] \dec1tree\ \\ \end{array}\bigg)
=\psi( | \vee_x | \vee_y |) = \overline{\lc\psi(|)\rc} \vee_x
\overline{\lc \psi(|)\rc} \vee_y \overline{\lc \psi(|)\rc} = x
\otimes y.$$
%
%$$ \psi(\arbreBCd\ \ \ )=\psi( | \vee_x | \vee_y |)
%= \overline{\lc\psi(|)\rc} \vee_x \overline{\lc \psi(|)\rc}
%\vee_y \overline{\lc \psi(|)\rc} = x \otimes y.$$
%

We now describe a submodule of $\ncshao(T(V))$ to be identified with the
image of $\psi$ in Theorem~\mref{thm:freetri}.

\begin{defn} {\rm Let $X=M(\Omega)$.
A $y\in \frakX_\infty$ is called a
{\bf dendriform triword (TW)} if it satisfies the following {\em
additional} properties.
\begin{enumerate}
%\item There is no subword $\lc \bfone \rc$ in the word;
\item $y$ is not in $\lc \frakX_\infty \rc$;
\item There is no subword $\lc \lc \frakx \rc \rc$ with $\frakx\in \frakX_\infty$
in the word;
%\item There is no subword of the form
%$\frakx_1\lc \frakx_2 \rc \ox_2$ with $\frakx_1,\ox_2\in X$ and
%$\frakx_2\in \frakX_\infty$.
\end{enumerate}
We let $TW(V)$ be the subspace of $\ncshao(T(V))$ generated by the
dendriform triwords. }
%The set of dendriform diwords is denoted by $DW(V)$.
\end{defn}

For example
$$ x_0\lc x_1\lc x_2 \rc \rc,
    \lc x_0\rc x_1\lc x_2 \rc,  \lc x_0\rc x_1\lc x_2 \rc x_3 \lc x_4\rc,
    x_0\otimes x_1$$
are dendriform triwords while
$$
%\lc \bfone \rc,
\lc \lc x_1\rc \rc,  %x_1\lc x_2 \rc x_3,
    \lc x_1 \lc x_2 \rc x_3\rc  $$
are in $\frakX_\infty$ but not dendriform triwords.

Equivalently, TWs can be characterized in terms of the
decomposition (\mref{eq:words3}). For subsets $Y,Z$ of
$\frakX_\infty$, define \allowdisplaybreaks{
\begin{eqnarray}
S(Y,Z)&=&\Big( \bigcup_{r\geq 1} (Y\lc Z\rc)^r \Big) \bigcup
    \Big(\bigcup_{r\geq 0} (Y\lc Z\rc)^r  Y\Big) \notag \\
&&    \bigcup \Big( \bigcup_{r\geq 1} \lc Z\rc (Y\lc Z\rc)^r \Big)
    \bigcup \Big( \bigcup_{r\geq 0} \lc Z\rc (Y\lc Z\rc)^r Y\Big).
\mlabel{eq:twords}
\end{eqnarray}}
Then define $S_0(V)=M(X)$. For $n\geq 0$, inductively define
\begin{equation}
 S_{n+1}(V)=S(M(X),S_n(V)).
%= (M(X)\lc T_n(V) \rc ) \bigcup (\lc S_n(V) \rc M(X))
%    \bigcup \lc S_n(V)\rc M(X) \lc S_n(V)\rc.
\mlabel{eq:tris}
\end{equation}
Then  $S_\infty: = \cup_{n\geq 0} S_n(V)$ is the set of dendriform triwords and
$TW(V)=\oplus_{\frakx\in S_\infty} \bfk \frakx.$

Theorem~\mref{thm:tri} follows from the following theorem.
\begin{theorem}
\begin{enumerate}
\item $\psi: \DT(V)\to \ncshao(T(V))$ is a homomorphism of
dendriform trialgebras. \mlabel{it:homt}
\item $\psi=T(j)$, the
morphism of dendriform trialgebras induced by $j:V\to
\ncshao(T(V))$. \mlabel{it:agreet}
\item $\psi(\DT)=DT(V)$.
\mlabel{it:imaget}
\item $\psi$ is injective. \mlabel{it:injt}
\end{enumerate}
\mlabel{thm:freetri}
\end{theorem}
\begin{proof} The proof is similar to Theorem~\mref{thm:freedend}.
For the lack of a uniform approach for both cases, we give
some details.

\mref{it:homt}
we first note that the operations $\prec$ and $\succ$ can be equivalently defined
as follows without using $|\prec T$, etc.
Let $T\in T_{i,X}, U\in T_{j,X}$ with $i\geq 1, j\geq 1$. Then
$T=T^{(0)}\vee_{x_1}\cdots \vee_{x_m} T^{(m)}$ and
$U=U^{(0)}\vee_{y_1}\cdots \vee_{y_n} U^{(n)}$, define
\begin{eqnarray*}
T\prec U &=& \left \{ \begin{array}{ll}
T^{(0)}\vee_{x_1} \cdots \vee_{x_m} (T^{(m)} \star U), & {\rm\ if\ } T^{(m)}\neq |,\\
T^{(0)}\vee_{x_1} \cdots \vee_{x_m} U, & {\rm\ if\ } T^{(m)}= |
\end{array} \right . \\
T\succ U &=& \left \{ \begin{array}{ll}
(T \star U^{(0)}) \vee_{y_1} \cdots \vee_{y_n} U^{(n)},& {\rm if\ } U^{(0)}\neq |,\\
T  \vee_{y_1} \cdots \vee_{y_n} U^{(n)},& {\rm if\ } U^{(0)}= |
\end{array} \right . \\
T\spr U &=& \left \{ \begin{array}{ll}
T^{(0)}\vee_{x_1} \cdots \vee_{x_m} (T^{(m)} \star U^{(0)})
\vee_{y_1} \cdots \vee_{y_n} U^{(n)},& {\rm if\ } T^{(m)}\neq |, U^{(0)}\neq |,\\
T^{(0)}\vee_{x_1} \cdots \vee_{x_m} U^{(0)}
\vee_{y_1} \cdots \vee_{y_n} U^{(n)},& {\rm if\ } T^{(m)}= |, U^{(0)}\neq |,\\
T^{(0)}\vee_{x_1} \cdots \vee_{x_m} T^{(m)}
\vee_{y_1} \cdots \vee_{y_n} U^{(n)},& {\rm if\ } T^{(m)}\neq |, U^{(0)}= |
%T^{(0)}\vee_{x_1} \cdots \vee_{x_m} |
%\vee_{y_1} \cdots \vee_{y_n} U^{(n)},& {\rm if\ } T^{(m)}\neq |, U^{(0)}= |,\\
\end{array} \right .
\end{eqnarray*}

Now we use induction on $i+j$ to prove
\begin{eqnarray}
&&\psi(T\prec U) = \psi(T) \prec_{R} \psi(U), \
\psi(T\succ U) = \psi(T) \succ_{R} \psi(U), \ \\
&&\psi(T\spr U) = \psi(T) \spr_{R} \psi(U).
\mlabel{eq:homt}
\end{eqnarray}
Here $R:=R_{T(V)}$ is the Rota-Baxter operator on $\ncshao(T(V))$.
Since $i+j\geq 2$, we can first take $i+j=2$. Then
$T=|\vee_x |$, $U=|\vee_y |$. So by Eq. (\mref{eq:psi}),
$$ \psi(T\prec U) = \psi( (|\vee_x |) \prec U)
= \psi( | \vee_x U) = x \lc \psi(U)\rc = x \lc y \rc
=x \prec_{R} y.$$
We similarly have $\psi(T\succ U)= x\succ_{R} y$ and
$$ \psi(T\spr U)=\psi((|\vee_x |) \spr (| \vee_y |))
= \psi( |\vee_x | \vee_ y |) = x\ot y = x \spr_{R} y.$$

Assume Equations (\mref{eq:homt}) hold for $T\in T_{i,X},\ U\in T_{j,X}$ with
$i+j\geq k\geq 2$. Then we also have
\begin{eqnarray}
 \psi(T\star U)&=&\psi (T\prec U+T\succ U +T\spr U) \notag \\
&=&    \psi(T) \prec_{R} \psi(U)+
    \psi(T) \succ_{R} \psi(U)+
    \psi(T) \spr_{R} \psi(U) \mlabel{eq:start}\\
&=& \psi(T)\star_{R} \psi(U). \notag
\end{eqnarray}

Here $\star_{R}=\prec_{R} +\succ_{R}+\spr_{R}.$ Consider $T,U$
with $m+n=k+1$. We consider two cases of
$T=T^{(0)}\vee_{x_1}\cdots \vee_{x_m} T^{(m)}$. Since $U\neq |$,
we have $\overline{\lc T^{(m)} \star U\rc}=\lc T^{(m)} \star U\rc$
if $T^{(m)}\neq |$, and $\overline{\lc U\rc}=\lc U\rc$ if
$T^{(m)}= |$.

{\bf Case 1.} If $T^{(m)}\neq |$, then
\begin{eqnarray*}
\psi(T\prec U) &=& \psi (T^{(0)}\vee_{x_1} \cdots
\vee_{x_m} (T^{(m)}\star U))
\ \ {\rm (definition\ of\ } \prec {\rm )} \\
&=& \overline{\lc \psi(T^{(0)}) \rc} x_1 \cdots x_m \lc \psi(T^{(m)} \star U)\rc
\ \ {\rm (definition\ of\ } \psi {\rm )} \\
&=& \overline{\lc \psi(T^{(0)}) \rc} x_1 \cdots x_m
    \lc \psi(T^{(m)}) \star_R \psi(U)\rc
\ \ {\rm (induction\ hypothesis\ (\mref{eq:start}))} \\
&=& \overline{\lc \psi(T^{(0)}) \rc} x_1 \cdots x_m
    \lc \psi(T^{(m)})\rc \lc \psi(U)\rc
\ \ {\rm (relation~(\mref{eq:RB}))} \\
&=& \psi (T^{(0)}\vee_{x_1} \cdots \vee_{x_m}
T^{(m)})\prec_R \psi(U)
\ \ {\rm (defintion\ of\ } \psi {\rm )}\\
&=& \psi(T)\prec_R \psi(U).
\end{eqnarray*}

{\bf Case 2.} If $T^{(m)}=|$, then
\begin{eqnarray*}
\psi(T\prec U) &=& \psi (T^{(0)}\vee_{x_1} \cdots
\vee_{x_m} U)
\ \ {\rm (definition\ of\ } \prec {\rm )} \\
&=& \overline{\lc \psi(T^{(0)}) \rc} x_1 \cdots x_m \lc \psi(U)\rc
\ \ {\rm (definition\ of\ } \psi {\rm )} \\
&=& \psi (T^{(0)}\vee_{x_1} \cdots \vee_{x_m} T^{(m)}) \lc \psi(U)\rc
\ \ {\rm (defintion\ of\ } \psi {\rm )}\\
&=& \psi(T)\prec_R \psi(U).
\end{eqnarray*}
This proves $\psi(T\prec U)=\psi(T)\prec_R \psi(U)$. We similarly prove
$\psi(T\succ U) =\psi(T)\succ_R \psi(U)$ and
$\psi(T\spr U) =\psi(T)\spr_R \psi(U)$.
Thus $\psi$ is a homomorphism in $\DT$.

\mref{it:agreet} follows from the uniqueness of the
morphism $\DT(V)\to \ncshao(T(V))$ of dendriform trialgebra extending the map
$i:V\to \ncshao(T(V))$.

\mref{it:imaget} We only need to prove $TW(V)\subseteq  \psi(\DT(V))$
and $\psi(\DT(V)) \subseteq TW(V)$. To prove the former, we prove
$S_n(V)\subseteq \psi(\DT(V))$ by induction on $n$.

When $n=0$, $S_n(V)=X$ so the inclusion is clear.
Suppose the inclusion holds for $1\leq n\leq k$. Then by the definition of
$S_{k+1}(V)$ in Eq.~(\mref{eq:tris}), an element of $S_{k+1}(V)$ has length
greater or equal to 2. We apply induction on its length.
If the length is 2, then it is one of the following two cases.

i) It is $\frakx \lc \frakx'\rc$ with $\frakx\in X$, $\frakx'\in
S_{k}(V)$. Then  it is $\frakx \prec_R \frakx'$
which is in $\psi(\DT(V))$ by the induction hypothesis and the
consequence from part \mref{it:homt} that $\psi(\DT(V))$ is a
sub dendriform algebra.

ii) It is $\lc \frakx \rc \frakx'$ with
$\frakx\in S_{k}(V)$  and
$\frakx'\in X$. Then the same proof works.

Suppose all elements of $S_{k+1}$ with length $\leq q$ and $\geq 2$ are in
$\psi(\DT(V))$. Consider an element $\frakx$ of $S_{k+1}$ with length $q+1$.
Then $q+1\geq 3$. If $q+1=3$, we again have two cases.

i) $\frakx=\lc \ox_1\rc \frakx_2 \lc \ox_3\rc$ with $\ox_1,\ox_2\in S_n(V)$
and $\frakx_1\in X$. Then it is
$(\ox_1 \succ_R \frakx_2) \prec_R \ox_3.$
By induction hypothesis on $n$, $\ox_1$ and $\ox_3$ are in the sub dendriform
dialgebra $\psi(\DT(V))$. So the element itself is in $\psi(\DT(V))$.

ii) $\frakx=\frakx_1 \lc \ox_2 \rc \frakx_3$ with $\frakx_1,\frakx_3\in X$
and $\ox_2\in S_n(V)$. Then
$\frakx=\frakx_1 \cdot_R (\ox_2 \succ \frakx_3)$ which is in $\psi(\DT(V))$.

If $q+1\geq 4$, then $\frakx$ can be expressed as the concatenation of
$\frakx_1$ and $\frakx_2$ of lengths at least two and hence are in
$TW(V)$. By induction hypotheses, $\frakx_1$ and $\frakx_2$ are in $\psi(\DT(V))$.
Therefore $\frakx = \frakx_1 \spr_R \frakx_2$ is in $\psi(\DT(V))$.

This completes the proof of the first inclusion. The proof of the second
inclusion follows from a similar induction on the degree of trees in $\DT(V)$.

\mref{it:injt} By the definition of $\psi$ and part \mref{it:image}, $\psi$ gives
a one-one correspondence between
$\cup_{n\geq 0} T_{n,X}$ as a basis of $\DT(V)$ and $TW(V)$ as a basis of
$\psi(\DT(V))$. Therefore $\psi$ is injective.
\end{proof}

%\addcontentsline{toc}{section}{\numberline {}References}
%\bibliography{reference}

\begin{thebibliography}{abcdsfgh}

\mbibitem{Ag1} M. Aguiar, {Pre-Poison algebras},
               {\em Lett. Math. Phys.}, {\bf 54}, (2000), 263-277.

\mbibitem{Ag2} M. Aguiar, {Infinitesimal Hopf algebras},
         {\em Contemporary Mathematics}, {\bf{267}}, (2000), 1-29.

\mbibitem{Ag3} M. Aguiar, {On the associative analog of Lie bialgebras},
               {\em Journal of Algebra}, {\bf{244}}, (2001), 492-532.

\mbibitem{A-L} M. Aguiar and J.-L. Loday, {Quadri-algebras},
                   {\em{J. Pure Applied Algebra}}, {\bf 191}, (2004), 205-221,
                   arXiv:math.QA/03090171

\mbibitem{A-M} M. Aguiar and W. Moreira, {Combinatorics of the free Baxter algebra,}
    {\em Electronic Journal of Combinatorics} {\bf 13}(1) (2006), R17: 38 pp, arXiv:math.CO/0510169

\mbibitem{A-S1} M. Aguiar and F. Sottile, {Structure of the Loday-Ronco Hopf algebra of trees},
                   {\em J. Algebra} {\bf 295} (2006), 473-511.

\mbibitem{A-S2} M. Aguiar and F. Sottile, {Cocommutative Hopf algebras of permutations and trees},
                    {\em J. of Algebraic Combinatorics} {\bf 22} (2005), 451-470, arXiv:math.QA/0403101.

\mbibitem{A-G-K-O} G. E. Andrews, L. Guo, W. Keigher and K. Ono,
    Baxter algebras and Hopf algebras,
    {\em Trans. Amer. Math. Soc.,} {\bf 355} (2003), 4639-4656.

\mbibitem{At} F. V. Atkinson, {Some aspects of Baxter's functional equation},
   {\em J. Math. Anal. Appl.}, {\bf{7}}, (1963), 1-30.

\mbibitem{Ba} G. Baxter, {An analytic problem whose solution follows from a simple algebraic identity,}
   {\em Pacific J. Math.}, {\bf 10}, (1960), 731-742.

\mbibitem{B-D} A. A. Belavin and V. G. Drinfeld, {Solutions of the classical Yang-Baxter equation
      for simple Lie algebras},
      {\em Funct. Anal. Appl.}, {\bf 16}, (1982), 159-180.

\mbibitem{Ca} P. Cartier, {On the structure of free Baxter algebras},
         {\em Adv. in Math.}, {\bf{9}}, (1972), 253-265.

\mbibitem{Ch} F. Chapoton, { Un th\'{e}or\`{e}me de Cartier-Milnor-Moore-Quillen pour les big\`{e}bres
             dendriformes et les alg\`{e}bres braces},
             {\em J. Pure Appl. Alg.}, {\bf 168}, (2002), 1-18.

\mbibitem{C-K1} A. Connes and D. Kreimer, { Renormalization in quantum field theory and
             the Riemann-Hilbert problem. I. The Hopf algebra structure of graphs
             and the main theorem.},
             {\em Comm. Math. Phys.}, {\bf 210}, (2000), no. 1, 249-273.

\mbibitem{C-K2} A. Connes and D. Kreimer, { Renormalization in quantum field theory and
                        the Riemann-Hilbert problem. II. The $\beta$-function, diffeomorphisms and
                        the renormalization group.},
                        {\em Comm. Math. Phys.}, {\bf 216}, (2001), no. 1, 215-241.

\mbibitem{EF1} K. Ebrahimi-Fard, { Loday-type algebras and the Rota-Baxter relation,}
                        {\em Letters in Mathematical Physics}, {\bf 61}, no. 2, (2002), 139-147.

\mbibitem{EF2} E. Ebrahimi-Fard, { On the associative Nijenhuis relation},
                        {\em The Electronic Journal of Combinatorics}, Volume 11(1), R38, (2004).

\mbibitem{E-G-G-V} K. Ebrahimi-Fard, Jos\'e M. Gracia-Bond\'{\i}a,
L. Guo and J. C. V\'arilly, Combinatorics of renormalization as
matrix calculus, {\em Physics Letters B}, {\bf 19} (2006), 552-558, arXiv:hep-th/0508154.

\mbibitem{shuf} K. Ebrahimi-Fard and L. Guo, { Quasi-shuffles, Mixable Shuffles and Hopf Algebras},
              {\em J. Algebraic Combinatorics} {\bf 24} (2006), 83-101, arXiv:math.RA/0506418.

\mbibitem{prod} K. Ebrahimi-Fard and L. Guo, {On the products and dual of binary,
    quadratic, regular operads}, {\em J. Pure and Applied Algebra}, {\bf 200} (2005), 293-317,
                             arXiv:math.RA/0407162.

\mbibitem{unit} K. Ebrahimi-Fard and L. Guo, Coherent unit actions on
    operads and Hopf algebras, arXiv:math.RA/0503342.

\mbibitem{mat} K. Ebrahimi-Fard and L. Guo, On matrix
representation of renormalization in prerturbative quantum field theory, arXiv:hep-th/0508155

\mbibitem{mzv} K. Ebrahimi-Fard and L. Guo, Rota-Baxter algebras and
    multiple zeta values, arxiv:math.NT/0601558.

\mbibitem{free} K. Ebrahimi-Fard and L. Guo, Free Rota-Baxter algebras and rooted trees,
    arXiv: math.RA/0510266.

\bibitem{EGfields06}
    K.~Ebrahimi-Fard and L.~Guo,
    {Rota-Baxter Algebras in Renormalization of Perturbative Quantum Field Theory,}
    {\em Fields Inst. Commun.}, {\bf 50}, (2007),~47-105,
    arXiv:hep-th/0604116.

\mbibitem{E-G-K2} K. Ebrahimi-Fard, L. Guo and D. Kreimer, {Integrable Renormalization II:
           the General case}, {\em Annales Henri Poincare} {\bf 6} (2005), 369-395.

\mbibitem{E-G-K3} K. Ebrahimi-Fard, L. Guo and D. Kreimer, { Spitzer's Identity and the Algebraic
          Birkhoff Decomposition in pQFT},
         {\em J. Phys. A: Math. Gen.}, {\bf 37}, (2004), 11037-11052.

\bibitem{egm2006}
    K. Ebrahimi-Fard, L. Guo and D. Manchon,
    {Birkhoff type decompositions and the Baker-Campbell-Hausdorff recursion,}
    {\em Commun. Math. Phys.}, {\bf 267}, (2006), 821-845.
    arXiv:math-ph/0602004.

\bibitem{ek2005}
    K.~Ebrahimi-Fard and D.~Kreimer,
    Hopf algebra approach to Feynman diagram calculations,
    {\em J.~Phys.~A}, {\bf 38}, (2005), R385-R406.
    arXiv:hep-th/0510202.

\bibitem{em2006}
    K. Ebrahimi-Fard and D. Manchon,
    {On matrix differential equations in the Hopf algebra of renormalization,}
    {\em Adv.~Theor.~Math.~Phys.}, {\bf 10}, (2006), 879-913.
    arXiv:math-ph/0606039.

\bibitem{EMP07}
    K.~Ebrahimi-Fard, D.~Manchon and F.~Patras,
    {A Bohnenblust-Spitzer identity for noncommutative Rota-Baxter algebras
    solves Bogoliubov's counterterm recursion,}
    arXiv:0705.1265v1 [math.CO].


\bibitem{KDF2007}
    K.~Ebrahimi-Fard, D.~Manchon and F.~Patras,
    {New identities in dendriform algebras,}
    arXiv:0705.2636v1 [math.CO].


\mbibitem{Fo} L. Foissy, { Les alg\`{e}bres de Hopf des arbres enracin\'{e}s d\'{e}cor\'{e}s II},
                        {\em Bull. Sci. Math.}, {\bf 126}, (2002), 249-288.

\mbibitem{Fra1} A. Frabetti, { Dialgebra homology of associative algebras},
                {\em C. R. Acad. Sci. Paris}, {\bf 325}, (1997), 135-140.

\mbibitem{Fra2} A. Frabetti, {Leibniz homology of dialgebras of matrices},
                {\em J. Pure Appl. Alg.}, {\bf 129}, (1998), 123-141.

\mbibitem{Fr} B. Fresse, Koszul duality of operads and homology of partition
posets, in  Homotopy theory: relations with algebraic geometry, group
cohomology,  and algebraic $K$-theory,  {\em Contemp. Math.}, {\bf 346}
(2004) Amer. Math. Soc., Providence, RI, 115--215.

\mbibitem{Gu} L. Guo, {Baxter algebras and differential algebras},
                        in ``Differential algebra and related topics", (Newark, NJ, 2000),
                        World Sci. Publishing, River Edge, NJ, (2002), 281-305.

\mbibitem{Gu2} L. Guo,
    { Baxter algebras and the umbral calculus,}
    {\em Adv. in Appl. Math.,} {\bf 27} (2001), 405-426.

\mbibitem{Gu5} L. Guo, { Baxter algebras, Stirling numbers and partitions}, {\em J. Algebra Appl.}
 {\bf 4} (2005), 153-164, arXiv:math.AC/0402348.

\mbibitem{Gu6} L. Guo, Operated semigroups, Motzkin paths and rooted trees, preprint, 2007.

\mbibitem{G-K1} L. Guo, W. Keigher, {Baxter algebras and shuffle products},
                        {\em Adv. Math.}, {\bf 150}, (2000), 117-149.

\mbibitem{G-K2} L. Guo,  W. Keigher, { On free Baxter algebras:
    completions and the internal construction,}
    {\em Adv. Math.} {\bf 151} (2000), 101--127.

\mbibitem{GK3} L. Guo,  W. Keigher, On differential Rota-Baxter algebras, arXiv: math.RA/0703780

\mbibitem{G-S} L. Guo and W. Yu Sit, {Enumeration of Rota-Baxter words}, in: Proceedings ISSAC 2006, Genoa, Italy, ACM Press, arXiv: math.RA/0602449.

\mbibitem{G-Z} L. Guo and B. Zhang, Renormalization of multiple
zeta values, arXiv:math.NT/0606076.

\mbibitem{Hol1} R. Holtkamp, {Comparison of Hopf algebras on trees},
                        {\em Arch. Math.}, (Basel) {\bf 80}, (2003), 368-383.

\mbibitem{Hol2} R. Holtkamp, {On Hopf algebra structures over operads},
                        preprint, July 2004, arXiv:math.RA/0407074.

\mbibitem{Ka} C. Kassel, {Quantum Groups}, Springer-Verlag, 1995.

\mbibitem{Kr1} D. Kreimer, {On the Hopf algebra structure of perturbative quantum field theories},
                        {\em Adv. Theor. Math. Phys.}, {\bf{2}}, (1998), 303-334.

\mbibitem{Le1} P. Leroux, {Ennea-algebras}, {\em J. Algebra}, {\bf{281}}, (2004), 287-302.

\mbibitem{Le2} P. Leroux, {Construction of Nijenhuis operators and dendriform trialgebras}, {\em Int. J. Math. Math. Sci.} (2004), no. 40-52, 2595-2615, arXiv:math.QA/0311132

\mbibitem{Le3} P. Leroux, {On some remarkable operads constructed from Baxter operators},
         arXiv:math.QA/0311214)

\mbibitem{Lo1} J.-L. Loday, {Dialgebras},
                        in Dialgebras and related operads, {\em Lecture Notes in Math.},
                        {\bf{1763}}, (2001), 7-66.(preprint 2001, arXiv:math.QA/0102053)

\mbibitem{Lo2} J.-L. Loday, {Scindement d'associativit\'{e} et alg\`{e}bres de Hopf},
                       {\em Proceedings of the Conference in honor of Jean Leray}, Nantes (2002),
                       S\'{e}minaire et Congr\`{e}s (SMF), {\bf{9}}, (2004), 155-172.

\mbibitem{Lo3} J.-L. Loday, {Arithmetree},
                        {\em J. Algebra}, {\bf258}, (2002), 275-309.

\mbibitem{Lo4} J.-L. Loday, {Completing the operadic butterfly}, {\em Georgian Math J.} {\bf 13} (2006), 741-749.


\mbibitem{L-R1} J.-L. Loday and M. Ronco, {Trialgebras and families of polytopes,}
                     in ``Homotopy Theory: Relations with Algebraic Geometry, Group
                     Cohomology, and Algebraic K-theory" Contemporary Mathematics, 346, (2004).
                     (preprint May 2002, arXiv:math.AT/0205043)

\mbibitem{L-R2} J.-L. Loday and M. Ronco, {Alg\`{e}bre de Hopf colibres},
                         {\em C. R. Acad. Sci. Paris}, {\bf 337}, (2003), 153-158.

\mbibitem{L-R3} J.-L. Loday and M. Ronco, {On the structure of cofree
    Hopf algebras,} {\em J. reine angew. Math.} {\bf 592} (2006) 123-155.

\mbibitem{MP1}
    D. Manchon and S. Paycha,
    {Shuffle relations for regularized integrals of symbols,} {\em Comm. Math. Phys.} {\bf 270} (2007), 13--51, arXiv:math-ph/0510067.

\mbibitem{MP2}
    D. Manchon and S. Paycha,
    {Chen sums of symbols and renormalized multiple zeta values,}
    arXiv:math.NT/0702135.


\mbibitem{Reu} C. Reutenauer, Free Lie Algebras, Oxford University
Press, 1993.

\mbibitem{R-STS1} A. G. Reyman and M. A. Semenov-Tian-Shansky,
    {Reduction of Hamitonian systems, affine Lie algebras and Lax
    equations. I}, {\em Invent. Math.}, {\bf 54}, (1979), 81-100.

\mbibitem{R-STS2} A. G. Reyman and M. A. Semenov-Tian-Shansky,
                    {Group theoretical methods in the theory of finite
                    dimensional integrable systems}, in: {\em Encyclopedia of mathematical
                    science, v.16: Dynamical Systems VII}, Springer, (1994), 116-220.

\mbibitem{Ron} M. Ronco, {Eulerian idempotents and Milnor-Moore theorem
                                  for certain non-cocommutative Hopf algebras,}
                                  {\em J. Algebra}, {\bf 254}, (2002), 152-172.

\mbibitem{Ro1} G. Rota, {Baxter algebras and combinatorial identities I,}
                                 {\em Bull. Amer. Math. Soc.,} {\bf 5}, 1969, 325-329.

\mbibitem{Ro2} G. Rota, {Baxter operators, an introduction,}
                                In: ``Gian-Carlo Rota on Combinatorics, Introductory papers
                                and commentaries", Joseph P.S. Kung, Editor,
                                Birkh\"{a}user, Boston, 1995.

\mbibitem{Ro3} G.-C. Rota and D. Smith, {Fluctuation theory and Baxter algebras},
  {\em Istituto Nazionale di Alta Matematica}, {\bf IX},
 (1972) 179-201. Reprinted in: ``Gian-Carlo Rota on Combinatorics:
 Introductory papers and commentaries'', J.P.S. Kung Ed.,
 Contemp. Mathematicians, Birkh\"auser Boston, Boston,
                                              MA, 1995.

\mbibitem{Sp} F. Spitzer, {A combinatorial lemma and its application to probability theory,}
                        {\em Trans. Amer. Math. Soc.}, {\bf 82}, (1956), 323-339.


\mbibitem{STS1} M. A. Semenov-Tian-Shansky, {What is a classical $r$-matrix?},
                                            {\em Funct. Ana. Appl.}, {\bf 17}, no.4., (1983), 259-272.

\mbibitem{zhao}
    J. Zhao,
    {Renormalization of multiple $q$-zeta values,}
    arXiv:math.NT/0612093.

\end{thebibliography}

\end{document}